\documentclass[11pt]{article}      
\usepackage{graphicx}
\usepackage{latexsym}
\usepackage{amssymb}
\usepackage[mathscr]{eucal}
\usepackage{color}
\newcommand{\Div}{\mbox{\rm div}\,}

\newcommand{\supp}{\mbox{\rm supp}\,}

\newcommand{\Int}[2]{{\displaystyle \int_{ #1}^{ #2}}}

\newcommand{\Lim}[1]{{\displaystyle \lim_{ #1}}}

\newcommand{\Sup}[1]{{\displaystyle \sup_{#1}}}

\newcommand{\Frac}[2]{\displaystyle{\frac{\displaystyle{#1}}{\displaystyle{#2}}}}

\newcommand{\beea}{\begin{eqnarray}}

\newcommand{\eeea}{\end{eqnarray}}

\newcommand{\ms}{\medskip\smallskip}

\newcommand{\bfe}{{\mbox{\boldmath $e$}} }

\newcommand{\bfz}{{\mbox{\boldmath $z$}} }

\newcommand{\0}{{\mbox{\boldmath $0$}} }

\newcommand{\BF}{\begin{footnotesize}}

\newcommand{\EF}{\end{footnotesize}}

\newcommand{\pde}[2]{{\displaystyle \frac{\mbox{$\partial #1$}}{\mbox{$\partial #2$}}}}

\newcommand{\ode}[2]{{\displaystyle \frac{\mbox{$d #1$}}{\mbox{$d #2$}}}}

\newcommand{\bi}{\begin{itemize}}

\newcommand{\ei}{\end{itemize}}

\newcommand{\ed}{\end{document}}

\newcommand{\be}{\begin{equation}}

\newcommand{\ba}{\begin{array}}

\newcommand{\ea}{\end{array}}

\newcommand{\ee}{\end{equation}}

\newcommand{\eeq}[1]{\label{eq:#1}\end{equation}}

\newcommand{\real}{{\mathbb R}}

\newcommand{\nat}{{\mathbb N}}

\newcommand{\bfpsi}{\mbox{\boldmath $\psi$}}

\newcommand{\bfchi}{\mbox{\boldmath $\chi$}}

\newcommand{\bfomega}{\mbox{\boldmath $\omega$}}

\newcommand{\bfxi}{\mbox{\boldmath $\xi$}}

\newcommand{\bfx}{\mbox{\boldmath $x$}}

\newcommand{\bfy}{\mbox{\boldmath $y$}}

\newcommand{\bfv}{{\mbox{\boldmath $v$}} }

\newcommand{\bfu}{{\mbox{\boldmath $u$}} }

\newcommand{\bfw}{{\mbox{\boldmath $w$}} }

\newcommand{\bff}{{\mbox{\boldmath $f$}} }

\newcommand{\bfa}{{\mbox{\boldmath $a$}} }

\newcommand{\bfc}{{\mbox{\boldmath $c$}} }

\newcommand{\bfA}{{\mbox{\boldmath $A$}} }

\newcommand{\bfQ}{{\mbox{\boldmath $Q$}} }

\newcommand{\bfI}{{\mbox{\boldmath $I$}} }

\newcommand{\bfh}{{\mbox{\boldmath $h$}} }

\newcommand{\cala}{{\cal A}}

\newcommand{\cald}{{\cal D}}

\newcommand{\calf}{{\cal F}}

\newcommand{\cali}{{\cal I}}

\newcommand{\cals}{{\cal S}}

\newcommand{\calt}{{\cal T}}

\newcommand{\bfV}{{\mbox{\boldmath $V$}} }

\newcommand{\bfU}{{\mbox{\boldmath $U$}} }

\newcommand{\bfF}{{\mbox{\boldmath $F$}} }

\newcommand{\bfb}{{\mbox{\boldmath $b$}} }

\newcommand{\bfg}{{\mbox{\boldmath $g$}} }

\newcommand{\bfn}{{\mbox{\boldmath $n$}} }

\newcommand{\half}{\mbox{$\frac{1}{2}$}}

\def\Bbb R{\real}

\def\hat{\widehat}

\def\tilde{\widetilde}

\def\bar{\overline}

\newcommand{\bfcalb}{\mbox{\boldmath ${\cal B}$}}

\newcommand{\bfcalg}{\mbox{\boldmath ${\cal G}$}}

\newcommand{\bfcalh}{\mbox{\boldmath ${\cal H}$}}

\newcommand{\bfcalf}{\mbox{\boldmath ${\cal F}$}}

\newcommand{\ED}{\end{description}}

\def\tag{\renewcommand{\theequation}}

\newcommand{\Footnote}{~\footnote}

\newcommand{\Br}{\begin{rem}\begin{rm}}

\newcommand{\Er}{\end{rm}\end{remark}}

\newtheorem{lemm}{Lemma}[section]

\newtheorem{theo}{Theorem}[section]

\newtheorem{rem}{Remark}[section]

\newtheorem{coro}{Corollary}[section]

\newtheorem{exe}{\footnotesize{Exercise}}[section]

\newcommand{\Be}{\begin{exe}\begin{footnotesize}\begin{rm}}

\newcommand{\EE}[1]{\end{rm}\end{footnotesize}\label{exe:#1}\end{exe}}

\newcommand{\Bt}{\begin{theo}\begin{sl}}
\newtheorem{proposition}{Proposition}[section]
\newcommand{\Bp}{\begin{proposition}\begin{sl}}
\newcommand{\EP}[1]{\end{sl}\label{proposition:#1}\end{proposition}}
\newcommand{\propref}[1]{{\rm Proposition \ref{proposition:#1}}}

\newcommand{\Et}{\end{sl}\end{theorem}}

\newcommand{\Bl}{\begin{lemm}\begin{sl}}

\newcommand{\El}{\end{sl}\end{lemma}}

\newcommand{\eqref}[1]{{\rm (\ref{eq:#1})}}

\newcommand{\Bc}{\begin{coro}\begin{sl}}

\newcommand{\Ec}{\end{sl}\end{coro}}

\newcommand{\ET}[1]{\end{sl}\label{theo:#1}\end{theo}}

\newcommand{\EL}[1]{\end{sl}\label{lemm:#1}\end{lemm}}

\newcommand{\theoref}[1]{{\rm Theorem \ref{theo:#1}}}

\newcommand{\ER}[1]{\end{rm}\label{rem:#1}\end{rem}}

\newcommand{\EC}[1]{\end{sl}\label{coro:#1}\end{coro}}

\newcommand{\lemmref}[1]{{\rm Lemma \ref{lemm:#1}}}

\newcommand{\essup}[1]{{\rm ess}\,{{\displaystyle \sup_{\hspace*{-6mm}{#1}}}}\,}

\begin{document}
\title{Navier-Stokes Flow past a Rigid Body\\ that Moves by Time-Periodic Motion}
\author{Giovanni P Galdi \smallskip \\ { \small Department of Mechanical Engineering and Materials Science}\\ { \small University of Pittsburgh, USA}}
\date{}
\maketitle
\begin{abstract}
We study existence,  uniqueness  and asymptotic spatial behavior of  time-periodic strong solutions to the Navier-Stokes equations in the exterior of a rigid body, $\mathscr B$, moving  by time-periodic motion of given period $T$, when the data are sufficiently regular and small. Our contribution improves  all previous ones in several directions. For example,  we allow both translational, $\bfxi$, and angular, $\bfomega$,  velocities of $\mathscr B$ to depend on time, and do not impose any restriction on the period $T$ nor on the averaged velocity, $\bar{\bfxi}$, of $\mathscr B$. If $\bfxi\not\equiv\0$ we  assume that $\bfxi$ and $\bfomega$ are both parallel to a  constant  direction, while no further assumption is needed if $\bfxi\equiv\0$. We also furnish the spatial asymptotic behavior of the velocity field, $\bfu$, associated to  such solutions. In particular, if $\mathscr B$ has a  net motion characterized by $\bar{\bfxi}\neq\0$, we then show that, at large distances from $\mathscr B$, $\bfu$ manifests a wake-like behavior in the direction $-\bar{\bfxi}$, entirely similar  to that of  the velocity field of the  steady-state flow occurring  when $\mathscr B$ moves with  velocity $\bar{\bfxi}$.

\end{abstract}
\renewcommand{\theequation}{{1}.\arabic{equation}}
\section{Introduction} The mathematical analysis of  time-periodic viscous flow around a {\em moving} body is a relatively new area of research.\footnote{For the case when the body is at rest, we refer the reader to \cite{GaSo,KMT,KoNa,Ma,MaPa,Y}.} The  basic problem at the foundation of this study can be  described as follows \cite{GaKy}. A rigid body, $\mathscr B$, moves  in an otherwise stagnant Navier-Stokes liquid, $\mathscr L$, that occupies the whole space outside $\mathscr B$. When referred to a  frame, $\calf$, attached to $\mathscr B$, the motion of $\mathscr B$ is time-periodic  of period $T$, and a given body force $\bfb$ of the same period $T$ may  be acting on $\mathscr L$. Then, the question to address is whether the corresponding flow of $\mathscr L$, in the frame $\calf$, will also be time-periodic and of period $T$ ({\em $T$-periodic}, in short). From the mathematical viewpoint, this means to find $T$-periodic solutions $(\bfu,p)$ to the following set of equations:
\be\ba{cc}\smallskip\left.\ba{ll}\medskip
\pde\bfu t-\bfV\cdot\nabla\bfu+\bfomega\times\bfu+\bfu\cdot\nabla\bfu=\Delta\bfu-\nabla {p}+\bfb\\
\Div\bfu=0\ea\right\}\ \ \mbox{in $\Omega\times (-\infty,\infty)$}\,,\\
\bfu=\bfV\,,\ \ \mbox{at\ $\partial\Omega\times (-\infty,\infty)$}\,,\ \ \Lim{|x|\to\infty}\bfu(x,t)=\0\,,\ \ t\in(-\infty,\infty)\,.
\ea
\eeq{0.1}
Here, $\bfu$ and $p$ represent velocity and pressure fields of $\mathscr L$,   and $\Omega$ is the spatial region outside  $\mathscr B$ and entirely occupied by $\mathscr L$. Furthermore, $\bfV=\bfV(t):=\bfxi(t)+\bfomega(t)\times\bfx$, where $\bfxi$ and $\bfomega$ are prescribed $T$-periodic functions denoting velocity of the center of mass and angular velocity of $\mathscr B$, respectively. Finally, $\bfb=\bfb(x,t)$ indicates a $T$-periodic body force acting on $\mathscr L$.
\par 
In spite of the conspicuous literature dedicated to the resolution of this problem, there is still a number of basic questions that, to date, are still unresolved, as we are going to present next. Well-posedness of \eqref{0.1} in the class of $T$-periodic solutions was first investigated in \cite{GS1} where it was shown that, if the data possess only a mild degree of regularity, there is at least one corresponding weak solution in the sense of Leray-Hopf. For more regular data, but of restricted ``size," in \cite{GS1} it is also constructed a solution that is strong in the sense of Ladyzhenskaya. The noteworthy aspect of these results is that, besides some regularity (and ``smallness" wherever it applies), {\em no further assumption} is made on the characteristic vectors $\bfxi$ and $\bfomega$. However, the unpleasant drawback is that, even in the class of strong solutions,  the {\em uniqueness} property  is not assured. This is due to the circumstance that  solutions constructed in \cite{GS1} carry very little information about their behavior as $|x|\to\infty$. In fact \eqref{0.1}$_4$ is satisfied only in a generalized sense, namely:\footnote{For notation, see the beginning of next section.}  
\be
\sup_{[0,T]}\|\bfu(t)\|_6<\infty\,.
\eeq{asbe}
The principal, but not only, difficulty in obtaining more relevant information for large $|x|$ is due to the presence in \eqref{0.1}$_1$ of the term $\bfomega(t)\times\bfx\cdot\nabla\bfu$ whose coefficient becomes unbounded as $|x|\to\infty$. For this reason, it is expected that the case  $\bfomega(t)\not\equiv\0$ will be the more challenging one. 
\par
More recently, the general problem of existence and uniqueness was  investigated and successfully settled by several authors  \cite{Eiter,EitKye1,EitKye,EitKye2,Ga,GaKy,GaKy1,GH,HG,MH,NTH,K1,Ngu1,Ngu},  with approaches and in classes of solutions different from those of \cite{GS1}. However, unlike \cite{GS1}, all results  there established hold under  hypotheses on $\bfxi$ and $\bfomega$ that are rather more restrictive  than just their regularity and smallness. Precisely, in \cite{EitKye,EitKye1,Ga,GaKy,GaKy1,HG,MH,K1} one assumes $\bfomega\equiv\0$,  $\bar{\bfxi}\neq\0$ (where the bar denotes  average over a period), and $\sup_t|\bfxi(t)|$, $\sup_t|\bfxi(t)-\bar{\bfxi}|$ ``small" enough. Still keeping $\bfomega\equiv\0$, in \cite{GARMA,GaLN} the assumption on $\bfxi$ was relaxed by requiring $\bfxi$ to be  only regular enough and ``small."  In \cite{Eiter,EitKye2} the request $\bfomega\equiv\0$ is  weakened  to $\bfomega={\bf const.}$, but on condition that $\bar{\bfxi}\neq\0$ and the period of $\bfxi$ and  $\bfb$ is  equal to $2\pi\kappa/|\bfomega|$, for some  $\kappa\in \mathbb Q\backslash\{0\}$. Moreover, at all times, $\bfxi(t)$ must be parallel to the (constant) direction of $\bfomega$,  with $\sup_t|\bfxi(t)|$, $\sup_t|\bfxi(t)-\bar{\bfxi}|$ and $|\bfomega|$  ``small" enough. The method there followed is based on maximal regularity properties  for  time-periodic problems with zero average, an approach  introduced in \cite{Ga} and successively fully developed and generalized  in \cite{K2,KS}. In \cite{GH,HG,MH,NTH,Ngu1,Ngu}, a different line of investigation  was undertaken  that relies  upon sharp ``$L^p-L^q$ estimates" of the relevant evolution operator in suitable spaces. However, this approach   requires   $\bfxi$ and $\bfomega$ to be {\em constant in time}, parallel and ``small." Nevertheless, unlike \cite{Eiter},  no restriction is imposed on the period $T$. 
\par
The main goal of this paper is to prove existence, {\em  uniqueness} and asymptotic spatial behavior of $T$-periodic solutions to \eqref{0.1}, under assumptions on $\bfxi$ and $\bfomega$ that are much more general than those requested in all papers cited above.\footnote{We are only interested to the case $\bfomega(t)\not\equiv\0$ because, otherwise, the problem has already been solved in \cite{GARMA,GaLN}.} Precisely, $\bfxi$ and $\bfomega$ are given $T$-periodic functions that, besides some  regularity and ``smallness", satisfy the following condition: 
\tag{H}
\be
\mbox{If $\bfxi(t)\not\equiv\0$, then both $\bfxi(t)$ and $\bfomega(t)$ are parallel to the same constant direction}
\eeq{H} 
It is worth emphasizing that if, instead, $\bfxi(t)\equiv\0$, {\em no further assumption is imposed on} $\bfomega(t)$.\footnote{An early version of this result with $\bfomega=\textbf{const.}$ was given in the preprint \cite{Garot}.} The significant feature of these conditions is that, unlike all previous contributions, both $\bfxi$ and $\bfomega$ are allowed to depend on time, and no restriction is imposed on the period $T$ nor on $\bar{\bfxi}$ being not zero. In particular, if the motion of $\mathscr B$ is only rotatory, its angular velocity can be an arbitrary (sufficiently regular and ``small")  periodic function of time. Under the above assumptions, we prove existence and uniqueness of $T$-periodic solutions in a class, $\mathscr C$, of functions belonging to suitable homogeneous Sobolev spaces (as  in \cite{GS1}),  with the further, critical property that the velocity field $\bfu(x,t)$ decays as $|x|^{-1}$, uniformly in time and even faster, outside a paraboloidal region (the ``wake"), if $\bar{\bfxi}\neq \0$. More precisely, in this case, $\bfu$ exhibits a wake-like behavior in the direction $-\bar{\bfxi}$, entirely analogous  to that of  the velocity field of the  steady-state flow that takes place when $\mathscr B$ moves with  velocity $\bar{\bfxi}$.
\par
Our method of proof stems from \cite{GARMA}, and does not require  advanced tools like maximal regularity or $L^p-L^q$ estimates. Instead, it  relies just on a skillful combination of the classical Galerkin method with uniform (spatial) estimates for solutions to  a suitable Oseen-like Cauchy problem, and  can be summarized as follows. The basic idea is to use a standard perturbative argument around solutions to the linear problem, $\mathscr{LP}$ (say), obtained from \eqref{0.1}  by suppressing the term $\bfu\cdot\nabla\bfu$ in \eqref{0.1}$_1$. The success of this approach is ensured provided we show existence,  uniqueness and continuous data-dependence for solutions to $\mathscr{LP}$ in the class $\mathscr C$. The proof of these properties is obtained in two steps.
In the first one, under suitable summability requirements on the data, we show the existence of a corresponding  solution in homogeneous Sobolev spaces. Successively, assuming in addition that the data decay pointwise in space   at a suitable rate and  uniformly in time, we prove that the above solution must enjoy a similar property as well. 
\par
To carry out  the first step we employ the well-known approach  that combines Galerkin method with the ``invading domains" technique; see, e.g., \cite{GS1}. Though this strategy is rather classical for the study of problems in unbounded domains, \cite[\S 17]{Ler} \cite[p. 180]{Lad}, \cite{Hey,GaSiRo}, in the case at hand its implementation presents a basic difficulty. Actually, in order to prove the desired asymptotic spatial behavior (objective of the second step) we need to construct a solution satisfying, in particular, the {\em crucial} requirement \renewcommand{\theequation}{{1}.\arabic{equation}}\setcounter{equation}{2}
\be
\essup{t}\int_{\Omega_R}\left|\pde\bfu t(t)\right|^q\le C,\ \ \mbox{ for some $q\ge 2$},  
\eeq{0.2}
where $\Omega_R=\Omega\cap\{|x|<R\}$, for some sufficiently large $R>0$, and $C$ is a positive constant depending on the data. 
When $\bfomega\equiv\0$, the proof of \eqref{0.2} with $q=2$   is rather standard by the above method. However, already when $\bfomega=\textbf{const.}\neq\0$, to show \eqref{0.2} is not trivial \cite{Garot}, due to the presence of the  term $\bfomega\times\bfx\cdot\nabla\bfu$, whose coefficient {\em grows  unbounded} at large spatial distances.  If $\bfomega$ is a generic (sufficiently smooth) $T$-paeriodic function, the situation is further complicated by the fact that the coefficient is also {\em time dependent}.    Thus, in such a case, by no means can the validity of \eqref{0.2}    be taken for granted,  let alone considered  an easy generalization of previous results. Nevertheless, by introducing new and non-trivial appropriate estimates we will be able to show \eqref{0.2} with $q=6$; see \propref{2.1}. 
\par
The solutions constructed in the first step  do not possess much  information as $|x|\to\infty$,  other than that given in \eqref{asbe}. As expcted, this is not enough to ensure the success of the perturbative argument.    Therefore,  in the second step, we investigate their asymptotic spatial behavior under the further condition that the data decay ``sufficiently fast" at large distances. This is done by reducing  $\mathscr{LP}$ to an analogous problem, $\mathscr{LP}_0$ (say) with $\Omega\equiv\real^3$ via a standard ``cut-off" procedure; see \eqref{1}. Successively, the treacherous term $\bfomega(t)\times\bfx\cdot\nabla\bfu$ is eliminated by means of a suitable time-dependent change of coordinates (see \eqref{2}) that brings $\mathscr{LP}_0$ into a Cauchy problem for an Oseen-like system of equations; see \eqref{5}. This is exactly the point (and the only one) where we need the assumption \eqref{H} on $\bfxi$ and $\bfomega$.  Classical results \cite[Theorem VIII.4.4]{Gab} then ensure existence of a solution, $(\bfv,{\sf p})$, to the Cauchy problem where $\bfv$ decays (at least) like $|x|^{-1}$ uniformly in time, on the assumption that the data decay ``sufficiently fast" and uniformly in time as well. The property \eqref{0.2} plays a fundamental role to guarantee the validity of such an assumption. Then, by uniqueness, also $\bfu$ must decay like $|x|^{-1}$, uniformly in time, and this  completes the second step; see \propref{2.2}.
\par 
Combining the findings of both steps, we then secure the desired results of existence, uniqueness and continuous data-dependence of solutions to  $\mathscr{LP}$ in the class $\mathscr C$; see \theoref{1.1}.  Thanks to this result and the functional properties of $\mathscr C$, we may finally use a standard contraction mapping argument that assures the existence and uniqueness of a solution to the full nonlinear problem \eqref{0.1} in $\mathscr C$, at least for  sufficiently smooth and ``small" data, and under the assumption \eqref{H} on $\bfxi$ and $\bfomega$; see \theoref{3.1}.  
\par
Before concluding this introductory section, some comments are in order as  whether \eqref{H}  is indeed necessary to prove well-posedness of  problem  \eqref{0.1} in the class of $T$-periodic solutions. I believe that if, as in the present paper, the class includes the pointwise asymptotic spatial behavior  which, on physical ground, is a very significant property, then  that hypothesis cannot be removed. However, if the class is enlarged to replace the pointwise behavior  of $\bfu$ (as $|x|^{-1}$) with the requirement that $\bfu$ belongs to the weak-$L^3$ space (as in the pioneering work \cite{Y}), there could be a chance that \eqref{H} might not be needed. This guess is also supported by the ``$L^p-L^q$ estimates" in weak spaces of the evolution operator associated to \eqref{0.1},  recently proved by {\sc T. Hishida} in \cite{Hish1,Hish2}, a fundamental tool for the success of the approach introduced in \cite{Y}.  
  
\par
The plan of the paper is as follows. After collecting some preliminary results in Section 2, in Section 3 we show the existence of a solution to $\mathscr{LP}$ in homogeneous Sobolev spaces. In the subsequent Section 4, we show the asymptotic spatial properties of these solutions and also their uniqueness and continuous data-dependence. In the final Section 5, we use the combined results of Section 3 and 4 to show existence, uniqueness and asymptotic spatial behavior of solutions to \eqref{0.1} under the assumption of data of restricted magnitude and with $\bfxi$ and $\bfomega$ satisfying the request \eqref{H}.

\renewcommand{\theequation}{\arabic{section}.\arabic{equation}}\setcounter{equation}{0}
\section{Preliminaries}
We begin to recall some notation. $\Omega$ will always denote  the complement of the closure of a bounded  domain $\Omega_0\subset\mathbb R^3$,  of class $C^2$. With the origin of the coordinate system in the interior of $\Omega_0$, we set $B_R=\{x\in\real^3:|x|<R\}$,  and if $R\ge R_*:=2\,{\rm diam}\,(\Omega_0)$,  
$
\Omega_R=\Omega\cap B_R\,,\ \ \Omega^R=\Omega\cap\bar{B_R}$.
For $A\subseteq \real^3$  a domain, by $L^q (A)$,  $1\leq q \leq \infty,$  
$W^{m,q}({A}),$ $W_0^{m,q}(A)$, $m \geq 0,$  $(W^{0,q}\equiv W^{0,q}_0\equiv L^q$), we indicate usual Lebesgue and Sobolev spaces, with norms $\|.\|_{q,A}$ and $\|.\|_{m,q,A}$, respectively.\Footnote{We shall use the same font style to denote scalar, vector and tensor
function spaces.} By $P$ we indicate  the (Helmholtz--Weyl) projector from $L^2(A)$ onto its subspace, $H(A)$, of solenoidal (vector) function with vanishing normal component (in the sense of trace) at $\partial A$.   
  We define $( u,v)_{A}:=\int_{A}u\cdot v$. With  
$D^{m,2}(A)$ we indicate the homogeneous Sobolev space of (equivalence classes of) functions $u$ with seminorm 
$ 
\sum_{|k|=m}\|D^k u\|_{2,A}<\infty\,.
$ 
In all the above notation,  the subscript ``$A$" will be omitted, unless otherwise   specified.
We also set
$
u_t:=\partial u/\partial t$,  $\partial_k u:=\partial{u}/\partial{x_k}\,,$
 and $D^2u=\{\partial_k\partial_ju\}$,  the matrix of the second derivatives.  A function $u:A\times \real\mapsto \real^3$ is 
{\em $T$-periodic}, $T>0$, if $u(\cdot,t+T)=u(\cdot\,t)$, for a.a. $t\in \real$,
 and we set
$
{\bar u}:=\frac{1}{T}\int_{0}^{T}u(t){\rm d}t\,.
$
Let $B$ be a function space endowed with seminorm $\|\cdot\|_B$, $r=[1,\infty]$, and $T>0$. $L^r(0,T;B)$ is the class of functions
$u:(0,T)\rightarrow B$ such that 
$$
\|u\|_{L^r(B)}\equiv\left\{\ba{ll}\smallskip\big( \Int{0}{T}\|u(t)\|_B^r \big)^{\frac 1r}<\infty, \ \ \mbox{if 
$r\in [1,\infty)\,;$}\\   
\essup{t\in[0,T]}\,\|u(t)\|_B <\infty, \ \ \mbox{if $r=\infty.$}
\ea\right.
$$
We also define
$$
W^{m,r}(0,T;B)=\Big\{u\in L^{r}(0,T;B): \textcolor{black}{\partial_t^ku\in L^{r}(0,T;B), \, k=1,\ldots,m}\Big\}\,.
$$
We shall simply write $L^r(B)$ for $L^r(0,T;B)$, etc. unless otherwise stated. Moreover, if $B\equiv\real^d$, $d\ge1$, we set $L^r(0,T;B)=L^r(0,T)$, etc.\par 
For $A:=\Omega,\real^3$,  $m=1,2$, and $\lambda\ge 0$ we set
$$
[\!] f [\!]_{\infty,m,\lambda,A}:=\Sup{(x,t)\in A\times (0,\infty)}\,|(1+|x|)^{m}(1+2\lambda \,s(x))^m f(x,t)|\,.
$$
where 
\be s(x)=|x|+x_1\,,\ \   x\in\real^3\,,
\eeq{esse} 
and  the subscript $A$ will be omitted, unless necessary.
\par
Finally, we denote by $c$ or $C$ a generic positive constant whose specific value is irrelevant and may change  even in the same line. When we want to emphasize the dependence of $c$ on a quantity $\rho$, we shall write $c_\rho$ or $c(\rho)$, and similarly for $C$. \smallskip\par
We now collect some preliminary results whose proof can be found in the literature. 
We begin with the following one \cite[Theorem III.3.1 and Exercise III.3.7]{Gab}.
\Bl Let $\cala$ be a bounded Lipschitz domain in $\real^3$, and let $f\in L^2(\cala)$ with $\int_{\cala}f=0$. Then the problem
\be
\Div\bfz=f\,\ \mbox{in $\cala$}\,, \ \bfz\in W_0^{1,2}(\cala)\,,\ \ \|\bfz\|_{1,2}\le C_0\,\|f\|_{2}\,,
\eeq{Bog}
for some $C_0=C_0(\cala)>0$ has at least one solution. Moreover, if also $f\in W_0^{1,2}(\cala)$, then $\bfz\in W_0^{2,2}(\cala)$ and $\|\bfz\|_{2,2}\le C_0\,\|f\|_{1,2}$.
Finally, if $f=f(t)$ with $ f_t\in L^\infty(L^2(\cala))$, then we have in addition $\bfz_t\in L^\infty(W_0^{1,2}(\cala))$ and
$$
\|\bfz_t\|_{L^\infty(W^{1,2})}\le C_0\,\|f_t\|_{L^\infty(L^{2})}\,.
$$ 
\EL{1.1_1}
The next result is proved in  \cite[Lemma 2.2]{GS1}. 
\Bl
Let $\bfxi,\bfomega \in W^{\textcolor{black}{2},2}(0,T)$ be $T$-periodic.There  exists a solenoidal, $T$-periodic function $\tilde{\bfu} \in W^{\textcolor{black}{1},2}(W^{m,q}),$ $m\in \mathbb{N},$ $q\in[1,\infty],$ such that 
$$\ba{ll}\medskip
\tilde{\bfu}(x,t)=\bfxi(t)+\bfomega\times\bfx\,,\ (t,\bfx)\in[0,T]\times\partial\Omega\,,\\ \medskip
\tilde{\bfu}(x,t)=0\,,\ \mbox{for all $(x,t)\in \bar{\Omega^\rho}\times [0,T]$}\,,\\ 
\| \tilde{\bfu} \|_{W^{2,2}(W^{m,q})}\leq C\,\left(\|\bfxi\|_{W^{2,2}(0,T)}+\|\bfomega\|_{W^{2,2}(0,T)}\right)
\,,\ea
$$
where $C=C(\Omega,m,q)$.  Moreover, for any $\varepsilon>0$, the field $\tilde{
\bfu}$ can be chosen in such a way that
\be
\left|\int_{\Omega_R}\bfv\cdot\nabla\tilde{\bfu}(t)\cdot\bfv\right|\le \varepsilon \|\nabla\bfv\|_2^2\,,\ \ \mbox{for all $\bfv\in H(\Omega_R)\cap W_0^{1,2}(\Omega_R)$}\,.
\eeq{LH}
\EL{ext}
We also have the following. 
\Bl
Let $\bfu\in H(\Omega_{R})\cap W_0^{1,2}(\Omega_R)$, and $\bfa,\bfb\in\real^3$. Then, the following properties hold
\begin{itemize}
\item[{\rm (i)}] $\bfa\times\bfx\cdot\nabla\bfu=\0\ \ \mbox{at $\partial B_R$}$\,;
\item[{\rm (ii)}]$
\big(\bfa \times \bfu - (\bfb+\bfa \times \bfx) \cdot \nabla \bfu\big) \in H(\Omega_{R})\,.
$
\end{itemize}
\EL{bello}
{\em Proof.} The proof of (i) is given in \cite[Lemma 3]{GaSiRo}. There, it was also shown that
$\big(\bfa \times \bfu - \bfa \times \bfx \cdot \nabla \bfu\big) \in H(\Omega_{R})$, so that, in order to prove (ii),  we only have to prove
$$
\bfb\cdot\nabla\bfu\in H(\Omega_R)\,.
$$ 
By \cite[Lemma III.2.1]{Gab} this is equivalent to showing that
\be 
(\bfb\cdot\nabla\bfu,\nabla\psi)=0\,,\ \ \mbox{for all $\psi\in D^{1,2}(\Omega_R)$\,.}
\eeq{FNC}
Since $\bfu\in H(\Omega_R)\cap W^{1,2}_0(\Omega_R)$, by \cite[Theorem III.4.1]{Gab} there is a sequence $\{\bfu_n\}$ of solenoidal functions from $C^\infty_0(\Omega_R)$ such that $\nabla\bfu_n\to\nabla\bfu$ in $L^2(\Omega_R)$. Clearly, as is shown by a simple integration by parts,
$$
(\bfb\cdot\bfu_n,\nabla\psi)=0\,,
$$
so that, by Schwarz inequality,
$$
|(\bfb\cdot\nabla\bfu,\nabla\psi)|=(\bfb\cdot\nabla(\bfu-\bfu_n),\nabla\psi)|\le |\bfb|\|\nabla(\bfu-\bfu_n)\|_2\|\nabla\psi\|_2\,,
$$
which, in turn, by letting $n\to\infty$, entails \eqref{FNC}\,.
\par\hfill$\square$\par
We conclude by recalling the next lemma  that ensures  suitable existence and uniqueness properties for a  linear Cauchy problem \cite[Theorem {VIII.4.4}]{Gab}
\Bl
Let $\bfcalg=\bfcalg(x,t)$ be a second-order tensor field defined in $\mathbb{R}^3 \times (0,\infty)$ such that
$$
\Div \bfcalg\in L^2(0,{\sf T};L^2)\,,
$$
for all ${\sf T}>0$,
and let $\bfh\in L^{\infty}(0,\infty;L^q)$, $q\in (3,\infty)$, with spatial support contained in $B_{\sf r}$, some ${\sf r}>0$.
Then, the problem 
\be
\ba{ll}\ms\left.\ba{rl}\ms
{\bfw}_{t}&=\Delta\bfw+\lambda\,\partial_1\bfw-\nabla{\sf q}+\Div\bfcalg+\bfh\\
\nabla\cdot\bfw&=0\ea\right\}\ \ \mbox{in $\real^3\times (0,\infty)$}\\
\bfw(x,0)=\0\,,
\ea
\eeq{VIII.4.1}
has one and only one solution such that
\begin{equation}
\bfw\in L^2(0,{\sf T};W^{2,2})\,, \ \bfw_t\in L^2(0,{\sf T};L^2)\,;\ \ \nabla{\sf q}\in L^{2}(0,{\sf T};L^2).
\label{estimunn1}
\end{equation}
If, in addition, $[\!]\bfcalg[\!]_{\infty,2,\lambda}<\infty$, then 
$$
[\!] \bfw [\!]_{\infty,1,\lambda}
+ \essup{t\geq 0}\| {\sf q}(t) \|_r
<\infty\,,
$$
for arbitrary $r\in (\frac{3}{2},\infty)$, 
and, setting 
$$
\cald:=[\!]\bfcalg[\!]_{\infty,2,\lambda} +\essup{t\geq 0}\|\bfh(t)\|_q\,,
$$
the following inequalities hold:
\begin{equation}
[\!] \bfw [\!]_{\infty,1,\lambda}\le C\,\cald\,,\ \  
 \essup{t\geq 0}\| {\sf q}(t) \|_r 
\leq C_1\,\cald\,,  
\label{estimunn2}
\end{equation}
with $C=C(q,{\sf r},\lambda_0)$, $C_1(q,r,{\sf r},\lambda_0)$ whenever ${\lambda}\in [0,\lambda_0]$, for some $\lambda_0>0$. 
\EL{VIII.4.4}
\setcounter{equation}{0}
\section{Linear Problem: Existence}
The main objective of this section and the following one is to prove existence, uniqueness and corresponding spatial asymptotic behavior of $T$-periodic solutions, in appropriate function classes, to the following linear problem
\be\ba{cc}\smallskip\left.\ba{ll}\medskip
\bfu_t-\bfV(t)\cdot\nabla\bfu+\bfomega(t)\times\bfu=\Delta\bfu-\nabla {p}+\bff\\
\Div\bfu=0\ea\right\}\ \ \mbox{in $\Omega\times (0,T)$}\\
\bfu(x,t)=\bfV(t)\,,\ \ (x,t)\in \partial\Omega\times [0,T]\,,
\ea
\eeq{4.7}
where $\bfV(t):=\bfxi(t)+\bfomega(t)\times\bfx$, and $\bff=\bff(x,t)$,  $\bfxi=\bfxi(t)$ and  $\bfomega=\bfomega(t)$, suitably prescribed $T$-periodic functions.
\par
Throughout, {\em we shall consider only the case $\bfomega(t)\not\equiv\0$}, since, otherwise, the problem has already been solved, even at the full nonlinear level, in  \cite{GARMA,GaLN}.  
\smallskip\par 
Our plan develops into two steps and goes as follows. In the first one, considered in this section, we  show, under suitable summability requirements on $\bff$, the existence of a corresponding  solution in Lebesgue and  homogeneous Sobolev spaces. Successively, in the next section, assuming in addition that $\bff$ decays pointwise in space   at a suitable rate  uniformly in time, we shall prove that the above solution must enjoy a similar property as well. We shall also show that, in this class, solutions are unique and depend continuously upon the data.
\par
To accomplish the first step, let $\cals=\{R_m, \ m\in\nat\}$ be an increasing,  unbounded sequence of positive numbers with 
$R_1 >\rho$ (see \lemmref{ext}), and let  $\{\Omega_{R}, \ R\in\cals\}$ be  the  sequence of bounded domains such that $\cup_{R\in \cals}\Omega_{R}=\Omega$.  
For each $R\in\cals$, we look for a $T-$periodic solution $(\bfu_R,p_R)$ to the problem
\be\ba{cc}\smallskip\left.\ba{ll}\medskip
(\bfu_R)_t-\bfV(t)\cdot\nabla\bfu_R+\bfomega(t)\times\bfu_R=\Delta\bfu_R-\nabla {p}_R+\bff\\
\Div\bfu_R=0\ea\right\}\ \ \mbox{in $\Omega_R\times (0,T)$;}\\
\bfu_R(x,t)=\bfV(t)\,,\ \ (x,t)\in \partial\Omega\times [0,T]\,;\ \ \bfu_R(x,t)=\0\,,\ \ (x,t)\in \partial B_R\times [0,T]\,.
\ea
\eeq{1.2_R}
Setting
$\bfv_R:=\bfu_R-\tilde{\bfu}$ with $\tilde{\bfu}$  given by \lemmref{ext}, we can equivalently rewrite \eqref{1.2_R} in the following form
\be\ba{cc}\smallskip\left.\ba{ll}\medskip
(\bfv_R)_t-\bfV(t)\cdot\nabla\bfv_R+\bfomega(t)\times\bfv_R=\Delta\bfv_R-\nabla {p}_R\\ \medskip\hspace*{5.1cm}-\tilde{\bfu}\cdot\nabla\bfv_R-{\bfv}_R\cdot\nabla\tilde{\bfu}+\tilde{\bff}\\
\Div\bfv_R=0\ea\right\}\ \ \mbox{in $\Omega_R\times (0,T)$}\\
\bfv_R(x,t)=\0\,,\ \ (x,t)\in \partial\Omega_R\times [0,T]\,,
\ea
\eeq{1.2}
where 
\be
\tilde{\bff}=\bff +  \Delta \tilde{\bfu} -  \tilde{\bfu}_t + \bfV \cdot \nabla \tilde{\bfu} - \bfomega \times \tilde{\bfu}:=\bff+\bff_c\,.
\eeq{sciu}
To show the existence of solutions to \eqref{1.2}--\eqref{sciu}, we shall employ Galerkin method with the orthonormal base of  $H(\Omega_R)$,  $\{\bfw _{Ri}\}_{i\in \mathbb{N}}$, 
constituted by the eigenfunctions of the Stokes operator:
\be
P\Delta \bfw_{Rj} = - \lambda_{Rj} \bfw_{Rj}\,,\ \ \ \bfw_{Rj}\in H(\Omega_R)\cap W_0^{1,2}(\Omega_R)\cap W^{2,2}(\Omega_R)\,.
\label{stokes}
\ee
We thus search for an ``approximating" solution to \eqref{1.2} of the form  
\be
\bfv_{Rk}(x,t)= \sum _{i=1}^{k}c_{Rki}(t)\bfw_{Ri} (x)\,,
\eeq{2.4_0} 
where  the coefficients $\bfc_{Rk}=\{c_{Rk1},\cdots,c_{Rkk}\}$  solve the following system of 
equations   
\be\ba{ll}\medskip
((\bfv_{Rk})_t,\bfw_{Rj})=(\nabla \bfv_{Rk} , \nabla \bfw_{Rj})_{\Omega_R} - ( \bfomega \times  \bfv_{Rk} , \bfw_{Rj})_{\Omega_R} + (\bfV \cdot \nabla \bfv_{Rk}, \bfw_{Rj})_{\Omega_R}\\
\hspace*{8.5cm} +(\tilde{\bff} , \bfw_{Rj})_{\Omega_R},\ \ j=1,\ldots,k\,. 
\ea
\eeq{ivpa}
Our next goal is to establish a number of estimates for $\bfv_{Rk}$ [respectively, $\bfv_R$] with bounds that are independent of $k$ [respectively, $R$]. 
In what follows, we will denote by $V_0$ a fixed positive number such that
\be
\|\bfxi\|_{W^{2,2}(0,T)}+\|\bfomega\|_{W^{2,2}(0,T)}\le V_0\,.
\eeq{V_0}
\par 
We begin to recall the following result whose proof is given in  \cite[Lemmas 3.1, 3.2, 4.1 and 4.3]{GS1}.\footnote{Actually, in the more general nonlinear context.} 
\Bl Let $\bff=\Div\bfcalf\in L^2(L^2)$ with $\bfcalf\in L^2(L^2)$, $\bfxi,\bfomega\in W^{1,2}(0,T)$ be $T$-periodic. Then, for each $k\in\nat$ problem \eqref{ivpa} has at least one $T$-periodic solution $\bfc_{Rk}=\bfc_{Rk}(t)$. Moreover, the approximating solution $\bfv_{Rk}$ satisfies the following uniform estimates
\be\ba{ll}\medskip
\Sup{t\in[0,T]}\left(\|\bfv_{Rk}(t)\|_{6,\Omega_R}+\|\nabla\bfv_{Rk}(t)\|_{2,\Omega_R}\right)+\|D^2\bfv_{Rk}(t)\|_{L^2(L^2(\Omega_R))}\\
\hspace*{3.7cm}
\le C\, \left(\|\bff\|_{L^2(L^2)}+\|\bfcalf\|_{L^2(L^2)}+\|\bfxi\|_{W^{1,2}(0,T)}+\|\bfomega\|_{W^{1,2}(0,T)}\right)\ea
\eeq{Art}
with $C=C(\Omega, V_0,T)$.
\EL{Ar}
We shall next prove additional  uniform estimates (in $k$) for the approximating solution \eqref{2.4_0}. 
\Bl Let the assumptions of\,        \lemmref{Ar} hold and suppose, in addition, $\bff,\bfcalf\in W^{1,2}(L^2)$ and $\bfxi,\bfomega\in W^{2,2}(0,T)$. Then $\bfv_{Rk}$ obeys the following bound
\be\ba{ll}\medskip
\Sup{t\in[0,T]}\|(\bfv_{Rk})_t(t)\|_{1,2,\Omega_R}+\|(\bfv_{Rk})_{tt}\|_{L^2(L^2(\Omega_R))}+\|D^2(\bfv_{Rk})_t\|_{L^2(L^2(\Omega_R))}\\ \hspace*{4cm}
\le C\, \left(\|\bff\|_{W^{1,2}(L^2)}+\|\bfcalf\|_{W^{1,2}(L^2)}+V_0\right)\,,\ea
\eeq{Art1}
where $C=C(R,T,V_0)$.
\EL{Ar1} 
{\em Proof.} In what follows, all norms and scalar products are taken in $\Omega_R$ which, therefore, will be omitted as a subscript. Moreover, we set $\bfv\equiv \bfv_{Rk}$. By taking the time-derivative of both sides of \eqref{ivpa}, dot-multiplying the resulting equations by $\dot{ c}_{Rkj}$,  summing over $j$ from 1 to $k$, and then integrating by parts over $\Omega_R$ we get:
\be
\frac12\ode{}t\|\bfv_t\|_2^2
=\!-\|\nabla\bfv_t(t)\|_2^2\!+\!(\dot{\bfV}\cdot\nabla\bfv-\dot{\bfomega}\times\bfv-\tilde{\bfu}_t\cdot\nabla\bfv\!-\!\bfv_t\cdot\nabla\tilde{\bfu}\!-\!\bfv\cdot\nabla\tilde{\bfu}_t\!+\!\bff_c,\bfv_t)-(\bfcalf_t,\nabla\bfv_t)
\,.
\eeq{oo}
By \eqref{LH} we may choose $\tilde{\bfu}$ such that
\be
|(\bfv_t\cdot\nabla\tilde{\bfu},\bfv_t)|\le \half\|\nabla\bfv_t\|_2^2\,.
\eeq{o2}
Also, from classical embedding theorems and \lemmref{ext}, we have
\be
\|\tilde{\bfu}\|_{L^\infty(W^{1,2})}+\|\tilde{\bfu}_t\|_{L^\infty(W^{1,2})}\le c\,\|\tilde{\bfu}\|_{W^{2,2}(W^{3,2})}\le c_1\,V_0\,,
\eeq{ET}
and so it follows that
$$
|(\dot{\bfV}\cdot\nabla\bfv+\tilde{\bfu}_t\cdot\nabla\bfv+\!\bfv\cdot\nabla\tilde{\bfu}_t\!-\bff_c,\bfv_t)|\le c\,\| \tilde{\bfu}\|_{W^{2,2}(W^{3,2})}\left(\|\nabla\bfv\|_2+\|\bfv\|_{2,\Omega_R}\right)\|\bfv_t\|_{2,\Omega_R}\,,
$$
where, we recall, $\Omega_R$ contains the bounded (spatial) support of $\tilde{\bfu}$. Thus, with the help of \lemmref{ext} and Poincar\'e inequality, we infer
\be
|(\dot{\bfV}\cdot\nabla\bfv-\dot{\bfomega}\times\bfv-\tilde{\bfu}_t\cdot\nabla\bfv-\bfv\cdot\nabla\tilde{\bfu}_t\!+\bff_c,\bfv_t)|\le c\,\|\nabla\bfv\|_2\|\nabla\bfv_t\|_2\,,
\eeq{o0}
with $c=c(R,V_0,T)$.
Finally, by Schwarz inequality,
\be
|(\bfcalf_t,\nabla\bfv_t)|\le \|\bfcalf_t\|_2\,\|\nabla\bfv_t\|_2\,.
\eeq{o3}
Replacing in \eqref{oo} the estimates \eqref{o2}--\eqref{o3}, we find
$$
\ode{}t\|\bfv_t\|_2^2+\|\nabla\bfv_t\|_2^2\le 
c\,\big(V_0+\|\nabla\bfv\|_2+\|\bfcalf_t\|_2\big)\|\nabla\bfv_t\|_2\,,
$$
which, in view of \lemmref{Ar}, implies
$$
\ode{}t\|\bfv_t\|_2^2+\|\nabla\bfv_t\|_2^2\le 
c\,\big(V_0+\|\bff\|_{L^2(L^2)}+\|\bfcalf\|_{L^{2}(L^2)}+\|\bfcalf_t\|_2\big)\|\nabla\bfv_t\|_2\,,
$$
Integrating both sides of the latter from 0 to $T$ and using the $T$-periodicity of $\bfv$ we readily conclude
\be
\|\nabla\bfv_t\|_{L^2(L^2)} \le 
c\,\big(V_0+\|\bff\|_{L^2(L^2)}+\|\bfcalf\|_{W^{1,2}(L^2)}\big)\,,
\eeq{1.11}
where $c=c(R,V_0,T)$ is independent of $\bfv$. We next take the time derivative of both sides of \eqref{ivpa}, dot-multiply the resulting equations by $-\lambda_{Rj}\dot{ c}_{Rkj}$ and use (\ref{stokes}). We then sum over $j$ from 1 to $k$, and integrate by parts over $\Omega_R$ to get 
\be\ba{ll}\medskip
\ode{}t\|\nabla\bfv_t\|_2^2
=-\|P\Delta\bfv_t\|_2^2 +\big(\bfV\cdot\nabla\bfv_t +\dot{\bfV}\cdot\nabla\bfv-\bfomega\times\bfv_t, P\Delta\bfv_t\big)\\
\hspace*{2.2cm}-\big(\tilde{\bfu}_t\cdot\nabla\bfv+\tilde{\bfu}\cdot\nabla\bfv_t+\bfv_t\cdot\nabla\tilde{\bfu}-\bfv\cdot\nabla\tilde{\bfu}_t+\tilde{\bff}_t,P\Delta\bfv_t\big)\,.
\ea
\eeq{1.12}
By arguing as before, and using Cauchy-Schwarz inequality, we can show
\be\ba{ll}\medskip
|\big(\tilde{\bfu}_t\cdot\nabla\bfv+\tilde{\bfu}\cdot\nabla\bfv_t+\bfv_t\cdot\nabla\tilde{\bfu}-\bfv\cdot\nabla\tilde{\bfu}_t+\tilde{\bff}_t,P\Delta\bfv_t\big)|\\
\hspace*{3cm}\le 
c\,[(V_0+\|\nabla\bfv\|_2+\|\nabla\bfv_t\|_2+\|\bff_t\|_2)]^2+\mbox{$\frac14$}\|P\Delta\bfv_t\|_2^2\,,\ea
\eeq{1.14}
with $c=c(V_0,T)$. 
Also, by Poincar\'e and Cauchy-Schwarz inequalities, we get
\be
\big|\big(\bfV\cdot\nabla\bfv_t +\dot{\bfV}\cdot\nabla\bfv-\bfomega\times\bfv_t, P\Delta\bfv_t\big)\big|\le c\,(\|\nabla\bfv\|_2^2+\|\nabla\bfv_t\|^2_2)+\mbox{$\frac14$}\|P\Delta\bfv\|_2^2\,,
\eeq{1.15}
with $c=c(V_0,R,T)$.
If we replace  \eqref{1.14} and \eqref{1.15} into \eqref{1.12}    we deduce
$$
\ode{}t\|\nabla\bfv_t\|_2^2
+\mbox{$\frac14$}\|P\Delta\bfv_t\|_2^2\\
\le  
c\,\big(V_0+\|\nabla\bfv\|_2+\|\nabla\bfv_t\|_2+\|\bff_t\|_2\big)^2\,,
$$
which, in turn, with the help of
 \eqref{Art} and \eqref{1.11} furnishes
\be
\ode{}t\|\nabla\bfv_t\|_2^2
+\mbox{$\frac14$}\|P\Delta\bfv_t\|_2^2\\
\le  
c\,\big(V_0+\|\bff\|_{W^{1,2}(L^2)}+\|\bfcalf\|_{W^{1,2}(L^2)}\big)^2\,,
\eeq{1.20}
with $c=c(V_0,T)$. We now observe that, by \eqref{1.11}, there is at least one $\bar{t}\in (0,T)$ such that
$$
\|\nabla\bfv_t(\bar{t})\|_{L^2} \le 
c\,\big(V_0+\|\bff\|_{L^2(L^2)}+\|\bfcalf\|_{W^{1,2}(L^2)}\big)\,,
$$
and so, integrating \eqref{1.20} between $\bar{t}$ and arbitrary $t>\bar{t}$, and exploiting the $T$-periodicity of $\bfv$, we readily get 
$$
\sup_{t\in[0,T]}\|\nabla\bfv_t({t})\|_{2} +\int_0^T\|P\Delta\bfv_t\|_2^2\le 
c\,\big(V_0+\|\bff\|_{W^{1,2}(L^2)}+\|\bfcalf\|_{W^{1,2}(L^2)}\big)\,.
$$
In turn, the latter, combined with the Poincar\'e inequality and the well-known inequality $\|D^2\bfw_{Ri}\|_2\le c\,\|P\Delta \bfw_{Ri}\|_2$, with $c=c(\Omega,R)$ \cite[Lemma IV.6.1]{Gab}, furnishes
\be
\sup_{t\in[0,T]}\|\bfv_t({t})\|_{1,2} +\int_0^T\|D^2\bfv_t\|_2^2\le 
c\,\big(V_0+\|\bff\|_{W^{1,2}(L^2)}+\|\bfcalf\|_{W^{1,2}(L^2)}\big)\,.
\eeq{NB}
Finally, we take the time-derivative of both sides of \eqref{ivpa}, dot-multiply the resulting equations by $\ddot{ c}_{Rkj}$ and  sum over $j$ from 1 to $k$ to get 
$$\ba{rl}\medskip
\|\bfv_{tt}\|_2^2=(\dot{\bfV}\cdot\nabla\bfv+&\bfV\cdot\nabla\bfv_t-\dot{\bfomega}\times\bfv-\bfomega\times\bfv_t-\tilde{\bfu}_t\cdot\nabla\bfv\\
&-\bfv_t\cdot\nabla\tilde{\bfu}-\tilde{\bfu}\cdot\nabla\bfv_t-\bfv\cdot\nabla\tilde{\bfu}_t+P\Delta\bfv_t,\bfv_{tt})+(\tilde{\bff}_t,\bfv_{tt})\,.
\ea
$$
Using Schwarz inequality on the right-hand side of this relation along with \eqref{ET}, and taking into account \eqref{sciu}, we show
\be
\|\bfv_{tt}\|_2\le c\, \big(\|\bfv_t\|_{2,2}+\|\bfv\|_{1,2}+\|{\bff}_t\|_2\big)\,,
\eeq{pp}
with $c=c(R,V_0,T)$. The lemma then follows from \eqref{NB}, \eqref{pp} and \lemmref{Ar}.
\par\hfill$\square$\par
Employing the results established in the two previous lemmas, we can now  prove the following one that guarantees the existence of solutions to the problem \eqref{1.2}, for each $R\in\cals$.
\Bl Suppose $\bff$, $\bfxi$ and $\bfomega$ satisfy the assumptions of \lemmref{Ar1}. Then, for any $R\in\cals$, problem \eqref{1.2} has one  $T$-periodic solution $(\bfv_R,p_R)$ such that\footnote{Notice that, since $p_R$ is $T$-periodic, from \eqref{1.2}$_1$ it follows that $\bfv_t$ is also $T$-periodic.\label{foot:T}} 
$$
\bfv_R\in W^{1,2}(W^{2,2}(\Omega_R))\cap W^{2,2}(L^2(\Omega_R))\,,\ \ p_R\in W^{1,2}(W^{1,2}(\Omega_R))\,. 
$$
Moreover, there is a constant $C=C(V_0,\Omega,T)$ independent of $R$, such that
\be\ba{ll}\medskip
\|\bfv_{R}\|_{L^\infty(L^6(\Omega_R))}+\|\nabla\bfv_{R}\|_{L^\infty(L^2(\Omega_R))}+\|D^2\bfv_{R}\|_{L^2(L^2(\Omega_R))}+\|\nabla p_R\|_{L^2(L^2(\Omega_R))}\\
\hspace*{3.6cm}
\le C\, \left(\|\bff\|_{L^2(L^2)}+\|\bfcalf\|_{L^2(L^2)}+\|\bfxi\|_{W^{2,2}(0,T)}+\|\bfomega\|_{W^{2,2}(0,T)}\right)\ea
\eeq{Artu}
\EL{erre}
{\em Proof.} In view of \lemmref{Ar} and \lemmref{Ar1}, the sequence of ``approximating solutions," $\{\bfv_{Rk}\}$, is bounded in the class $W^{1,2}(W^{2,2})\cap W^{2,2}(L^2)$, uniformly in $k$. Therefore,  one may find a subsequence  $\{\bfv_{Rk^\prime}\}$ and a function 
\be
\bfv_R\in W^{1,2}(W^{2,2}(\Omega_R))\cap W^{2,2}(L^2(\Omega_R))
\eeq{vR} 
such that $\bfv_{Rk^\prime}\to \bfv_R$ as $k^\prime\to\infty$ in appropriate topology. By a classical argument \cite{GS1}, we then show that there is $p_R\in L^2(W^{1,2}(\Omega_R))$ such that $(\bfv_R,p_R)$ is a $T$-periodic solution to \eqref{1.2}.  However, from \eqref{1.2}$_1$ and \eqref{vR} it follows that $p_R\in W^{1,2}(W^{1,2}(\Omega_R))$. Clearly, by \lemmref{Ar}, $\bfv$ satisfies \eqref{Artu}. As a result, in order to complete the proof of the lemma, it remains to show the estimate for $p_R$ in \eqref{Artu}. By the Helmholtz decomposition, it follows that 
\be
\tilde{\bff}=P\tilde{\bff}+\nabla \tilde p\,, \ \ \ \ \tilde p\in L^2(D^{1,2}),
\label{hd}
\ee 
where
$$
\|\nabla\tilde{p}\|_{L^2(L^2)}\le c\,\|\tilde{\bff}\|_{L^2(L^2)}\,,
$$
and $c=c(\Omega)$. From \eqref{sciu} and the latter, we deduce
\be
\|\nabla\tilde{p}\|_{L^2(L^2)}\le c\,(\|{\bff}\|_{L^2(L^2)}+\|\bfxi\|_{W^{2,2}(0,T)}+\|\bfomega\|_{W^{2,2}(0,T)})\,,
\eeq{RM}
where $c=c(\Omega,V_0)$.
Recalling that, by \lemmref{bello},
$$
\left((\bfxi+\bfomega\times \bfx)\cdot \nabla \bfv_R-\bfomega\times \bfv_R\right)\in H(\Omega_R)
\,,
$$
from \eqref{1.2}$_1$ we show that ${\sf p}_R:= p_R-\tilde p|_{\Omega_R}$ satisfies the following Neumannn problem (in the sense of distributions)
\be
\ba{ll}\ms
\Delta {\sf p}_R=\Div \bfF\,, \ \ \mbox{in $\Omega_R$} 
\\
\pde{{\sf p}_R}{n}=\bfF\cdot n\,, \ \ \mbox{at $\partial \Omega_R$},
\ea\label{buono1}
\ee
where $\bfF:=- \bfv_R \cdot \nabla \tilde{\bfu} 
- \tilde{\bfu}\cdot \nabla \bfv_R+ \Delta \bfv_R$. Formally multiplying both sides of (\ref{buono1})$_1$ by ${\sf p}_R$, integrating by parts over $\Omega_R$ and using (\ref{buono1})$_2$ we deduce
$$
\|\nabla {\sf p}_R\|^2_{2,\Omega_R}=(\bfF,\nabla {\sf p}_R)_{\Omega_R}\,,
$$
which, in turn, by Schwarz inequality, entails
\be
\|\nabla {\sf p}_R\|^2_{2,\Omega_R}\le\|\bfF\|_2^2
\eeq{PR}
From \eqref{ET} we get
\be
\|\bfF\|_{L^2(L^2)}\le c\, (\|\bfv_R\|_{L^2(W^{1,2}(K))}+\|P\Delta\bfv_R\|_{L^2(L^2(\Omega_R))})
\eeq{FK}
with $c=c(V_0,T)$ and where $K$ is the bounded (spatial) support of $\tilde{\bfu}$. Collecting (\ref{hd}), \eqref{RM}, \eqref{PR},  \eqref{FK} and \eqref{Artu} for $\bfv_R$ we show
\be
\|\nabla p_R\|_{L^2(L^2(\Omega_R))}
\le C\, \left(\|\bff\|_{L^2(L^2)}+\|\bfcalf\|_{L^2(L^2)}+\|\bfxi\|_{W^{2,2}(0,T)}+\|\bfomega\|_{W^{2,2}(0,T)}\right)
. 
\eeq{nece}
Finally, since $\bfF$ is $T$-periodic, by the (well known) uniqueness property for the problem (\ref{buono1}) we infer that $P_R$, and hence $p_R$, is $T$-periodic as well. The proof of the lemma is thus completed.
\par\hfill$\square$\par
Our next objective is to show that the solution obtained in \lemmref{erre}, in addition to \eqref{Artu}, satisfies some further estimates, uniformly with respect to $R$. 
In this regard, we need some preparatory results.
\Bl Let $\bff$, $\bfxi$,  $\bfomega$  and $(\bfv_R,p_R)$ be as in \lemmref{erre}. Then
\be\ba{ll}\medskip
\|(\bfv_R)_t-\bfV(t)\cdot\nabla\bfv_R+\bfomega(t)\times\bfv_R\|_{L^2(L^2(\Omega_R))}\\
\hspace*{2.5cm}\le c_1 \left(\|\bff\|_{L^2(L^2)}+\|\bfcalf\|_{L^2(L^2)}+\|\bfxi\|_{W^{1,2}(0,T)}+\|\bfomega\|_{W^{1,2}(0,T)}\right)\,,
\ea
\eeq{cz1}
where $c_1=c_1(\Omega,V_0,T)$. Moreover,
\be
\|\tilde{\bfu}_t\cdot\nabla\bfv_R+\tilde{\bfu}\cdot\nabla(\bfv_R)_t+(\bfv_R)_t\cdot\nabla\tilde{\bfu}+\bfv_R\cdot\nabla\tilde{\bfu}_t\|_2\le c_2\,(\|\nabla\bfv_R\|_2+\|\nabla(\bfv_R)_t\|_2)
\eeq{cz2}
with $c_2=c_2(\Omega,V_0,T)$.
\EL{Chop}
{\em Proof.} From \eqref{1.2} we have\footnote{All spatial norms of $\bfv_R$ and $p_R$ are meant to be taken in $\Omega_R$.}
\be
\|(\bfv_R)_t-\bfV(t)\cdot\nabla\bfv_R+\bfomega(t)\times\bfv_R\|_{L^2(L^2)}=
\|\Delta\bfv_R-\nabla {p}_R-\tilde{\bfu}\cdot\nabla\bfv_R-{\bfv}_R\cdot\nabla\tilde{\bfu}+\tilde{\bff}\|_{L^2(L^2)}\,.
\eeq{V0}
Clearly, by \eqref{ET},  \eqref{sciu}, and \lemmref{ext}, we infer
\be\ba{ll}\medskip
\|\Delta\bfv_R-\nabla {p}_R-\tilde{\bfu}\cdot\nabla\bfv_R+\tilde{\bff}\|_{L^2(L^2)}\\
\hspace*{.7cm}\le c\,\big(\|\nabla\bfv_R\|_{L^2(W^{1,2})}+\|\nabla p_R\|_{L^2(L^2)}+\|\bff\|_{L^2(L^2)}+\|\bfxi\|_{W^{2,2}(0,T)}+\|\bfomega\|_{W^{2,2}(0,T)}\big)\,,
\ea
\eeq{V1}
with $c=c(V_0,\Omega,T)$. Moreover, recalling that the support of $\tilde{\bfu}$ is contained in $\Omega_\rho$, and taking into account \eqref{ET}, we obtain
\be
\|\bfv_R\cdot\nabla\tilde{\bfu}\|_{L^2(L^2)}\le cV_0\|\bfv_R\|_{L^2(L^2(\Omega_\rho))}\le c_1V_0 \|\bfv_R\|_{L^2(L^6)}
\eeq{V2}
with $c_1$ independent of $R$. The inequality \eqref{cz1} is then a consequence of combining \eqref{V0}--\eqref{V2} with \eqref{Artu}. Finally, \eqref{cz2} readily follows from \eqref{ET} and, since $\bfv_R\equiv \0$ at $\partial\Omega$,  the Poincar\'e inequality applied to $\bfv_R$ and $(\bfv_R)_t$ on the domain $\Omega_\rho$.
\par\hfill$\square$
\par
\Bl Let $\bfxi,\bfomega$ and $\bfv_R$ be as in \lemmref{Chop}. The following inequality holds\footnote{All scalar products and norms are taken on the domain $\Omega_R$.}
\be
|\big((\bfb+\bfa\times\bfx)\cdot\nabla\bfv_R-\bfa\times\bfv_R,P\Delta(\bfv_R)_t\big)|\le c\, V_0\|\nabla\bfv_R\|_{1,2}(\|\nabla(\bfv_R)_t\|_2+\|P\Delta(\bfv_R)_t\|_2)\,,
\eeq{V3}
with either $\bfb=\bfxi$ and $\bfa=\bfomega$, or $\bfb=\dot{\bfxi}$ and $\bfa=\dot{\bfomega}$, and $c=c(\Omega)$. Moreover, for any $\eta>0$ there is $c=c(\Omega,V_0,\eta)$ such that
\be  
|\big((\bfxi+\bfomega\times\bfx)\cdot\nabla(\bfv_R)_t-\bfomega\times(\bfv_R)_t,P\Delta(\bfv_R)_t\big)|\le c\, \|\nabla(\bfv_R)_t\|_{2}^2+\eta\|P\Delta(\bfv_R)_t\|_2^2\,.
\eeq{V4}
\EL{V1}
{\em Proof.} For simplicity, we set $\bfv\equiv\bfv_R$. Clearly, by Schwarz inequality,
\be
|(\bfb\cdot\nabla\bfv,P\Delta\bfv_t)|\le V_0\|\nabla\bfv\|_2\|P\Delta\bfv_t\|_2\,.
\eeq{V5}
Furthermore, by \lemmref{bello}(ii), 
$$
\cali:=\big(\bfa\times\bfx\cdot\nabla\bfv-\bfa\times\bfv,P\Delta\bfv_t\big)=
\big(\bfa\times\bfx\cdot\nabla\bfv-\bfa\times\bfv,\Delta\bfv_t\big)\,,
$$
so that, integrating by parts over $\Omega_R$ and taking into account  \lemmref{bello}(i) and that $\bfv$ vanishes at $\partial\Omega_R$, we show
\be
\cali=-\Int{\partial\Omega}{}\bfn\cdot\nabla\bfv_t\cdot(\bfa\times\bfx\cdot\nabla\bfv)+\Int{\Omega_R}{}\nabla(\bfa\times\bfx\cdot\nabla\bfv-\bfa\times\bfv):\nabla\bfv_t\,.
\eeq{Vi}
Also with the help of classical trace theorems, from this relation we infer
\be
|\cali|\le c\,V_0\,\|\nabla\bfv\|_{1,2}(\|\nabla\bfv_t\|_2+\|D^2\bfv_t\|_2)+\big|\Int{\Omega_R}{}\nabla(\bfa\times\bfx\cdot\nabla\bfv):\nabla\bfv_t\big|\,,
\eeq{V6}
with $c=c(\Omega)$. We now have
\be
\Int{\Omega_R}{}\nabla(\bfa\times\bfx\cdot\nabla\bfv):\nabla\bfv_t=\Int{\Omega_R}{}\varepsilon_{pqk}a_q\partial_pv_i\partial_kv_{ti}+\Int{\Omega_R}{}\varepsilon_{pqr}a_qx_r\partial_p\partial_kv_i\partial_kv_{ti}:=\cali_1+\cali_2\,,
\eeq{V7}
where $\varepsilon_{pqr}$ is the alternating symbol.
Obviously,
\be
|\cali_1|\le c\,V_0\,\|\nabla\bfv\|_2\|\nabla\bfv_t\|_2\,,
\eeq{V8}
with $c$ numerical constant. Moreover,
$$\ba{rl}\medskip
\cali_2=\half\Int{\Omega_R}{}\bfa\times\bfx\cdot\nabla(\partial_t|\nabla\bfv|^2)&\!\!\!=\half\Int{\partial\Omega}{}\bfa\times\bfx\cdot\bfn\,\partial_t|\nabla\bfv|^2+\half\Int{\partial B_r}{}\bfa\times\bfx\cdot\bfn\,\partial_t|\nabla\bfv|^2\\
&\!\!\!=\Int{\partial\Omega}{}\bfa\times\bfx\cdot\bfn (\nabla\bfv:\nabla\bfv_t)\,.
\ea
$$
Therefore, again by trace theorems
\be
|\cali_2|\le c\,V_0\,\|\nabla\bfv\|_{1,2}(\|\nabla\bfv_t\|_2+\|D^2\bfv_t\|_2)\,.
\eeq{V9}
We now recall Heywood inequality \cite[Lemma 1]{Hey}
\be
\|D^2\bfw\|_2\le c\,(\|P\Delta\bfw\|_2+\|\nabla\bfw\|_2)\,,\ \ \bfw\in H(\Omega_R)\cap W^{1,2}_0(\Omega_R)\cap W^{2,2}(\Omega_R)\,,
\eeq{JGH}
where the constant $c$ is independent of $R$. Thus, the proof of \eqref{V3}  becomes  a consequence of \eqref{V5}, \eqref{V6}--\eqref{JGH}. To prove \eqref{V2} we observe that by Cauchy-Schwarz inequality
\be
|(\bfxi\cdot\nabla\bfv_t,P\Delta\bfv_t)|\le c\,V_0\|\nabla\bfv_t\|_2^2+\eta\,\|P\Delta\bfv_t\|_2^2\,,
\eeq{V10}
with $c=c(V_0)$, 
and that by \lemmref{bello}(ii), 
\be
\calt:=\big(\bfomega\times\bfx\cdot\nabla\bfv_t-\bfomega\times\bfv_t,P\Delta\bfv_t\big)=
\big(\bfomega\times\bfx\cdot\nabla\bfv-\bfomega\times\bfv,\Delta\bfv_t\big)
\eeq{V11}
As a result, integrating by parts over $\Omega_R$ and proceeding exactly as in \eqref{Vi}, we show
\be
\calt=-\Int{\partial\Omega}{}\bfn\cdot\nabla\bfv_t\cdot(\bfomega\times\bfx\cdot\nabla\bfv_t)+\Int{\Omega_R}{}\nabla(\bfomega\times\bfx\cdot\nabla\bfv_t-\bfomega\times\bfv_t):\nabla\bfv_t\,.
\eeq{V12}
By using an argument entirely analogous to that employed in the proof of \eqref{V3}, we may prove
\be
\big|\Int{\Omega_R}{}\nabla(\bfomega\times\bfx\cdot\nabla\bfv_t-\bfomega\times\bfv_t):\nabla\bfv_t\big|\le\half\Big|\Int{\partial\Omega}{}\bfomega\times\bfx\cdot\bfn \,|\nabla\bfv_t|^2\Big|+c\,V_0\|\nabla\bfv_t\|_2^2\,,
\eeq{V13}
where $c$ is a numerical constant. Therefore, from \eqref{V12} and \eqref{V13}, we deduce
\be
|\calt|\le c\,V_0(\|\nabla\bfv_t\|_{2,\partial\Omega}^2+\|\nabla\bfv_t\|_2^2)\,,
\eeq{V14}
with $c=c(\partial\Omega)$. Again by trace theorems, we know that for any $\eta>0$ there exists $c=c(\eta,\Omega)$ such that \cite[Exercise II.4.1]{Gab}
$$
\|\nabla\bfv_t\|_{2,\partial\Omega}^2\le c\,\|\nabla\bfv_t\|_2^2+\eta\|D^2\bfv_t\|_2^2
$$
and \eqref{V4} follows from the latter, \eqref{V14} and \eqref{JGH}.
\par\hfill$\square$\par
Collecting the previous two lemmas, we can now prove the desired further uniform estimate.
\Bl Let $\bff$, $\bfxi$,  $\bfomega$  and $(\bfv_R,p_R)$ be as in \lemmref{erre}. Then, the following estimate holds
\be\ba{ll}\medskip \|(\bfv_{R})\|_{L^\infty(L^\infty(\Omega_R))}+
\|(\bfv_{R})_t\|_{L^\infty(L^6(\Omega_R))}+\|\nabla(\bfv_{R})_t\|_{L^\infty(L^2(\Omega_R))}+\|D^2(\bfv_{R})_t\|_{L^2(L^2(\Omega_R))}\\

\hspace*{.2cm}
+\|\nabla (p_R)_t\|_{L^2(L^2(\Omega_R))}\le C\, \left(\|\bff\|_{W^{1,2}(L^2)}+\|\bfcalf\|_{L^2(L^2)}+\|\bfxi\|_{W^{2,2}(0,T)}+\|\bfomega\|_{W^{2,2}(0,T)}\right)\,,\ea
\eeq{V15}
where $C=C(\Omega,V_0,T)$.
\EL{EST}
{\em Proof.} As usual, we shall omit the subscript $\Omega_R$ in the various scalar products and norms, and set $\bfv\equiv\bfv_R$. Moreover, we
let
$$\bfU:=\bfv_t-\bfV\cdot\nabla\bfv-\bfomega\times\bfv\,,\ \ E:=\|\bfU\|_2^2
$$ 
We  take the time derivative of both sides of \eqref{1.2}$_1$ and dot-multiply both sides of the resulting equation, called \eqref{1.2}$_t$, by $\bfU$.
Integrating by parts as necessary and observing that, by \lemmref{bello}(ii), $(\bfV\cdot\nabla\bfv-\bfomega\times\bfv,\nabla p_t)=0$, we show
$$
\half\ode Et+\|\nabla\bfv_t\|_2^2=(\Delta\bfv_t,\bfV\cdot\nabla\bfv-\bfomega\times\bfv)-(\tilde{\bfu}_t\cdot\nabla\bfv+\tilde{\bfu}\cdot\nabla\bfv_t+\bfv_t\cdot\nabla\tilde{\bfu}+\bfv\cdot\nabla\tilde{\bfu}_t-\tilde{\bff}_t,\bfU)\,.
$$
By virtue of \lemmref{bello}(ii), \lemmref{Chop} and \lemmref{V1} the previous relation furnishes
$$
\half\ode Et+\|\nabla\bfv_t\|_2^2\le c\,\big[\|\nabla\bfv\|_{1,2}(\|\nabla\bfv_t\|_2+\|P\Delta\bfv_t\|_2)+(\|\nabla\bfv\|_2+\|\nabla\bfv_t\|_2+\|\tilde{\bff}_t\|_2)E^{\frac12}\big]\,,
$$
where $c=c(\Omega,V_0,T)$. Using, suitably, Cauchy-Schwarz inequality on the right-hand side of the latter, we get
\be
\ode Et+\|\nabla\bfv_t\|_2^2\le \eta\|P\Delta\bfv_t\|_2^2+c_1\,(\|\nabla\bfv\|_{1,2}^2+E+\|\tilde{\bff}_t\|_2^2)\,,
\eeq{V16}
where $\eta>0$ is arbitrary and $c_1=c_1(\Omega,V_0,T,\eta)$. Next, we dot-multiply both sides of \eqref{1.2}$_t$ by $P\Delta\bfv_t$ and integrate by parts as necessary. In this way we show
$$\ba{rl}\medskip
\half\ode{}t\|\nabla\bfv_t\|_2^2+\|P\Delta\bfv_t\|_2^2=&\!\!\!-(\dot{\bfV}\cdot\nabla\bfv-\dot{\bfomega}\times\bfv,P\Delta\bfv_t)-({\bfV}\cdot\nabla\bfv_t-{\bfomega}\times\bfv_t,P\Delta\bfv_t)\\
&\!\!\!+(\tilde{\bfu}_t\cdot\nabla\bfv+\tilde{\bfu}\cdot\nabla\bfv_t+\bfv_t\cdot\nabla\tilde{\bfu}+\bfv\cdot\nabla\tilde{\bfu}_t-\tilde{\bff}_t,P\Delta\bfv_t)\,.
\ea
$$
Again by \lemmref{Chop} and \lemmref{V1}, from this equation we deduce
$$\ba{ll}\medskip
\half\ode{}t\|\nabla\bfv_t\|_2^2+\!\|P\Delta\bfv_t\|_2^2\le\! c\big[\|\nabla\bfv\|_{1,2}(\|\nabla\bfv_t\|_2+\|P\Delta\bfv_t\|_2)\\
 \hspace{2cm}+(\|\nabla\bfv\|_2+\|\nabla\bfv_t\|_2+\|\tilde{\bff}\|_2)\|P\Delta\bfv_t\|_2+\|\nabla\bfv_t\|_2^2\big]+\mbox{$\frac14$}\|P\Delta\bfv_t\|_2^2\,.
\ea
$$
where $c=c(\Omega,V_0,T)$. Thus, using Cauchy-Schwarz inequality, we readily obtain
\be
\ode{}t\|\nabla\bfv_t\|_2^2+\|P\Delta\bfv_t\|_2^2\le c_2\,(\|\nabla\bfv_t\|_2^2+\|\nabla\bfv\|_{1,2}^2+\|\tilde{\bff}\|_2^2)
\eeq{V17}
with $c_2=c_2(\Omega,V_0,T)$. We now multiply both sides of \eqref{V16} by $2c_2$, choose $\eta=1/4c_2$ and sum, side by side, the resulting  inequality
and \eqref{V17}. We thus get
\be
\ode{}t\left(\|\nabla\bfv_t\|_2^2+2c_2E\right)+c_2\|\nabla\bfv_t\|_2^2+\half\|P\Delta\bfv_t\|_2^2\le c_3\,(E+\|\nabla\bfv\|_{1,2}^2+\|\tilde{\bff}\|_{W^{1,2}(L^2))})\,.
\eeq{V18}
Next, observe that, by \eqref{Artu} and \lemmref{V1} it follows that
\be
\int_0^T (E(t)+\|\nabla\bfv(t)\|_{1,2}^2){\rm d}t\le c\,\big[\|{\bff}\|_{W^{1,2}(L^2))}+\|\bfcalf\|_{L^2(L^2)}+\|\bfxi\|_{W^{2,2}(0,T)}+\|\bfomega\|_{W^{2,2}(0,T)}\big)\,,
\eeq{V19}
and that, by  \lemmref{ext} and \eqref{sciu},
\be
\|\tilde{\bff}\|_{W^{1,2}(L^2))}\le c\,\left(\|{\bff}\|_{W^{1,2}(L^2))}+\|\bfxi\|_{W^2(0,T)}+\|\bfomega\|_{W^2(0,T)}\right)
\eeq{V20}
with $c=c(\Omega,T,V_0)$. Therefore, 
Integrating both sides of \eqref{V18} from 0 to $T$, then using the $T$-periodicity of $E$ and $\nabla\bfv_t$,\footnote{See footnote \ref{foot:T}.} along with \eqref{V19},  \eqref{V20}, and \eqref{JGH}
we find
\be
\|\nabla\bfv_t\|_{L^2(L^2)}^2+\|D^2\bfv_t\|_{L^2(L^2)}^2\le c\,\big[\|{\bff}\|_{W^{1,2}(L^2))}+\|\bfcalf\|_{L^2(L^2)}+\|\bfxi\|_{W^{2,2}(0,T)}+\|\bfomega\|_{W^{2,2}(0,T)}\big]\,.
\eeq{V21}
Further, from \eqref{V19} and \eqref{V21} it follows that there is $\bar{t}\in (0,T)$ such that
\be
E(\bar{t})+\|\nabla\bfv_t(\bar{t})\|_2^2\le c\,\big[\|{\bff}\|_{W^{1,2}(L^2))}+\|\bfcalf\|_{L^2(L^2)}+\|\bfxi\|_{W^{2,2}(0,T)}+\|\bfomega\|_{W^{2,2}(0,T)}\big]\,,
\eeq{JL}
so that, integrating \eqref{V18} from $\bar{t}$ to $2T$ and using the $T$-periodicity,  \eqref{V19},  \eqref{V20}, and \eqref{JL} we conclude, in particular,
\be
\sup_{t\in[0,T]}\|\nabla\bfv_t(t)\|_2^2\le c\,\big[\|{\bff}\|_{W^{1,2}(L^2))}+\|\bfcalf\|_{L^2(L^2)}+\|\bfxi\|_{W^{2,2}(0,T)}+\|\bfomega\|_{W^{2,2}(0,T)}\big]\,.
\eeq{V22}
If we now extend $\bfv_t$ to 0 outside $\Omega_R)$ and continue to denote by $\bfv_t$ the extension, we may use the Sobolev inequality to get 
\be
\|\bfv_t\|_6\le c_\Omega\,\|\nabla\bfv_t\|_2\,.
\eeq{2.57}
Furthermore, by well--known embedding theorems (e.g. \cite[Corollary 5.16]{Adams}) it follows that
$$
\|\bfv\|_\infty\le c\,(\|\bfv\|_6+\|\nabla\bfv\|_2+\|D^2\bfv\|_2)
$$
with $c$ depending only on the regularity of $\Omega$. Thus, since by \lemmref{erre} and \eqref{V21} $\bfv\in L^\infty(L^6\cap D^{1,2})\cap W^{1,2}(D^{2,2})$, and $W^{1,2}(D^{2,2})\subset L^\infty(D^{2,2})$, the last displayed inequality furnishes
$$
\|\bfv\|_{L^\infty(L^\infty)}\le c\,(\|\bfv\|_{L^{\infty}(L^6)}+\|\nabla\bfv\|_{L^{\infty}(L^2)}+\|D^2\bfv\|_{W^{1,2}(D^{2,2})})\,,
$$
where $c$ is depending only on the regularity of $\Omega$.
As a result, the desired estimate for $\bfv$ follows from the latter, \lemmref{erre}  and \eqref{V21}, \eqref{V22}, \eqref{2.57}. Concerning the stated estimate for $\nabla p_t$, we see that its proof can be carried out exactly in the same way as \eqref{nece} and will, therefore, omitted.
\par\hfill$\square$\par
With the help of \lemmref{erre} and \lemmref{EST} we are now able to prove a general existence result of $T$-periodic solutions to \eqref{4.7}. To this end, define the function spaces
\be\ba{ll}\medskip
\mathscr U:=\left\{ \mbox{$T$-periodic}\ \bfu: \, \bfu\in W^{1,\infty}(L^6\cap D^{1,2})\cap W^{1,2}(D^{2,2})\cap L^\infty(L^\infty)\,;\ \Div\bfu=0\right\}\\
\hat{\mathscr P}:=\left\{\mbox{$T$-periodic}\ p:\, p\in L^\infty(L^6)\cap W^{1,2}(D^{1,2})\right\}\,. 
\ea
\eeq{V23}
Clearly, both $\mathscr U$ and $\mathscr P$ become Banach spaces when endowed with the norms
$$\ba{ll}\medskip
\|\bfu\|_{\mathscr U}:=\|\bfu\|_{L^\infty(L^\infty)}\!+\!\|\bfu\|_{W^{1,\infty}(L^6)}\!+\!\|\nabla\bfu\|_{W^{1,\infty}(L^2)}+\|D^2\bfu\|_{W^{1,2}(L^2)}\,;\\ 
\|p\|_{\hat{\mathscr P}}:=\|p\|_{L^\infty(L^6)}+\|\nabla p\|_{W^{1,2}(L^2)}\,.
\ea
$$
\par
Our next result represents the  main finding of this section, and establishes the existence of $T$-periodic solutions to problem \eqref{4.7} in suitable Lebesgue and homogeneous Sobolev spaces.  
\Bp Suppose $\bff$,  $\bfxi$ and $\bfomega$ satisfy the assumptions of \lemmref{Ar1}. Then, problem  \eqref{4.7} has at least one  solution $(\bfu,p)\in \mathscr U\times\hat{\mathscr P}$ such that 
\be
\|\bfu\|_{\mathscr U}+\|p\|_{\hat{\mathscr P}}
\le C \left(\|\bff\|_{W^{1,2}(L^2)}+\|\bfcalf\|_{W^{1,2}(L^2)}+\|\bfxi\|_{W^{2,2}(0,T)}+\|\bfomega\|_{W^{2,2}(0,T)}\right)\,,
\eeq{V29}
with $C=C(\Omega, T,V_0)$. 
\EP{2.1}
{\em Proof.} Set $\bfu_R:=\bfv_R+\tilde{\bfu}$, with $(\bfv_R,p_R)$ given in \lemmref{erre}. Then,  $(\bfu_R,p_R)$ is a solution to \eqref{1.2_R}. Because of the properties of the extension $\tilde{\bfu}$, it is at once recognized that  both estimates \eqref{Artu} and \eqref{V15} hold with $\bfv_R\equiv\bfu_R$, namely,
\be
\|\bfu_R\|_{L^\infty(L^\infty(\Omega_R))}\!+\!\|\bfu_R\|_{W^{1,\infty}(L^6(\Omega_R))}\!+\!\|\nabla\bfu_R\|_{W^{1,\infty}(L^2(\Omega_R))}+\|D^2\bfu_R\|_{W^{1,2}(L^2(\Omega_R))}\le c\,\mathscr D\,,
\eeq{2.60}
where $c=c(\Omega,T,V_0)$ and
$$
\mathscr D:=\|\bff\|_{W^{1,2}(L^2)}+\|\bfcalf\|_{W^{1,2}(L^2)}+\|\bfxi\|_{W^{2,2}(0,T)}+\|\bfomega\|_{W^{2,2}(0,T)}
.
$$ 
We want to let $R\to\infty$ ($R\in\cals$) and show that $(\bfu_R,p_R)$ tends (in suitable topology) to a solution $(\bfu,p)\in \mathscr U\times\hat{\mathscr P}$ to \eqref{4.7}. This can be done by an argument similar to that given in  \cite[Section 3]{Ga}. Let $\chi=\chi(s)$, $s>0$, be a smooth, non-increasing real function such that $\chi(s)=1$ for $s\le 1/2$ and $\chi(s)=0$ for $s\ge1$. For $R_m\in\cals$,  set
$
\chi_m(x)=\chi\left(|x|/R_m\right). 
$ We thus have, for all $x\in\Omega$,
$$
\chi_m(x)=\left\{\ba{ll}\medskip 1 &\ \mbox{if $|x|\le \half R_m$}\,,\\
0 &\ \mbox{if $|x|\ge R_m$\,.}\ea\right.
$$
Notice that 
\be
|D^\alpha\chi_m(x)|\le C R_m^{-|\alpha|}\,, \ \ |\alpha|\ge 0\,;\ \ {\rm supp}\,(\chi)\subseteq \{\half R_m\le|x|\le R_m\}\,,
\eeq{Chi}
where $C$ is independent of $x$ and $m$. Let $\bfu_m:=\bfu_{R_m}$. We observe that $\bfu_m$ satisfies the following inequality
\be
\Int{\half R_m\le |x|\le R_m}{}|\bfu_m(x)|^2|x|^{-2}\le 4\,\|\nabla\bfu_m\|_{2,\Omega_{R_m}}^2\,.
\eeq{INE}
In fact, let's extend $\bfu_m$ to 0 outside the ball of radius $R_m$ and continue to denote by $\bfu_m$ such an extension. Then, clearly, $\bfu_m\in D^{1,2}(\Omega)$ and $\lim_{|x|\to\infty}|\bfu_m(x)|=0$. Thus, \eqref{INE} follows from \cite[Theorem II.6.1(i)]{Gab}. 
 We now set $\hat{\bfu}_m:=\chi_m\bfu_m$. From \eqref{2.60},  \eqref{Chi} and \eqref{INE} we readily deduce that $\hat{\bfu}_m\in\mathscr U$ and 
$$
\|\hat{\bfu}_m\|_{\mathscr U}\le c\,\mathscr D
$$
where, from now onward, $c$ denotes a constant independent of $m$. Therefore, there exists $\bfu\in\mathscr U$ such that (possibly, along a subsequence)
\be
\hat{\bfu}_m\to\bfu\,,\ \ \mbox{weakly in $\mathscr U$}\,,
\eeq{2.63}
and, in addition,
\be
\|{\bfu}\|_{\mathscr U}\le c\,\mathscr D\,.
\eeq{2.64}
We shall show that
\be
\int_0^\tau\int_{\Omega}\left({\bfu}_{t}-\bfV(t)\cdot\nabla{\bfu}+\bfomega(t)\times{\bfu}-\Delta{\bfu}-\bff\right)\cdot\bfpsi=0\,,
\eeq{2.67}
for all $\bfpsi\in C^\infty_0(\Omega)$ with $\Div\bfpsi=0$, and all $\tau\in [0,T]$. Actually,    let us dot-multiply by $\bfpsi$ both sides of \eqref{1.2_R}$_1$ where   $R\equiv R_m$ is so large that $\Omega_{\frac12 R_m}$ contains the support, $K$, of $\bfpsi$. By integrating over $\Omega_{R_m}\times (0,\tau)$, $\tau\in [0,T]$ the resulting equation, (also by parts, when necessary), we get 
$$
\int_0^\tau\int_{K}\left(\bfu_{mt}-\bfV(t)\cdot\nabla\bfu_m+\bfomega(t)\times\bfu_m-\Delta\bfu_m-\bff\right)\cdot\bfpsi=0\,.
$$
Since $\bfu_m(x,t)=\hat{\bfu}_m(x,t)$ for all $(x,t)\in K\times [0,T]$, this relation furnishes
\be
\int_0^\tau\int_{\Omega}\left(\hat{\bfu}_{mt}-\bfV(t)\cdot\nabla\hat{\bfu}_m+\bfomega(t)\times\hat{\bfu}_m-\Delta\hat{\bfu}_m-\bff\right)\cdot\bfpsi=0\,.
\eeq{2.66}
Thus, passing to the limit $m\to\infty$ in \eqref{2.66} and employing \eqref{2.63} we arrive at \eqref{2.67}.
Now, 
since $\bfpsi$ is arbitrary in its class, from well-known results and the fact that $\bfu\in\mathscr U$, we deduce that there is $p\in L^2(W^{1,2}(\Omega_R))$, arbitrary $R>R_*$, such that $(\bfu,p)$ satisfy \eqref{4.7}$_1$  a.e. in $\Omega\times[0,T]$. Moreover, it is not difficult to prove that $\bfu$ obeys  \eqref{4.7}$_3$ a.e. in $\partial\Omega\times[0,T]$. Finally, for any $\phi\in C_0^\infty(\Omega)$ we have for all sufficiently large $m$ 
$$
\int_0^\tau\int_{\Omega}\hat{\bfu}_m\cdot\nabla\phi=\int_0^\tau\int_{\Omega_{R_m}}{\bfu}_m\cdot\nabla\phi=0\,,
$$ 
which, by \eqref{2.63}, implies that $\bfu$ satisfies \eqref{4.7}$_2$ a.e. in $\Omega\times[0,T]$. Thus, in order to complete the proof of the proposition, it remains to show the stated properties for the pressure field $p$. To this end, we recall the Hardy-type inequality \cite[Theorem II.6.1(i)]{Gab}
\be 
\|u \,|x|^{-1}\|_2\le c\|\nabla u\|_2\,,\ \ u\in L^6(\Omega)\cap D^{1,2}(\Omega)\,.
\eeq{2.68}
As a result, multiplying both sides of \eqref{4.7}$_1$ by $|x|^{-1}$ and bearing in mind that $\bfu\in \mathscr U$, we deduce a.e. in $[0,T]$
\be
|x|^{-1}\nabla p\in L^2(\Omega)\,.
\eeq{2.69}
Observing that
\be\ba{ll}\medskip
\Div(\bfV\cdot\nabla\bfu+\bfomega\times\bfu)=\Div\bfu_t=\Div\Delta\bfu=0\\
\bfV\cdot\nabla\bfu=0\,,\ \ \mbox{at $\partial\Omega$}\,;\ \ \Delta(\dot{\bfV}+\bfomega\times\bfV)=\0\,,
\ea
\eeq{purchia}
from \eqref{4.7}$_1$ we obtain for a.a. $t\in [0,T]$ and in distributional sense
\be
\Delta p=\Div\bff\ \ \mbox{in $\Omega$}\,;\ \ \pde{p}{n}=-\bfn\cdot\Delta\bfu\ \ \mbox{at $\partial\Omega$}\,.
\eeq{2.70} 
Classical results on the Neumann problem ensure the existence of at least one solution, $\hat{p}$,  to \eqref{2.70} such that
\be
\|\nabla\hat{p}\|_2+\|\nabla\hat{p}_t\|_2\le c\,\big(\|\bff\|_2+\|\bff_t\|_2+\|\bfcalf\|_2+\|\bfcalf_t\|_2+\|\Delta\bfu\|_2+\|\Delta\bfu_t\|_2\big)\,.
\eeq{2.71}
Setting ${\sf p}:=p-\hat{p}$, it follows that, for a.a. $t\in [0,T]$,
\be
\Delta {\sf p}=0\ \ \mbox{in $\Omega$}\,;\ \ \pde{\sf p}{n}=0\ \ \mbox{at $\partial\Omega$}\,.
\eeq{2.72}
Since $\hat{p}\in D^{1,2}(\Omega)$ and $p$ satisfies \eqref{2.69}, again by classical results on the Neumann problem, we have  \cite[Exercise V.3.6]{Gab}
\be
D^\alpha{\sf p}=O(|x|^{-|\alpha|-1})\,,\ \  |\alpha|\ge 0, \ \mbox{as $|x|\to\infty$}\,.
\eeq{2.73} 
Therefore, multiplying both sides of \eqref{2.72} by ${\sf p}$, integrating by parts over $\Omega_R$, and then letting $R\to\infty$, with the help of \eqref{2.73} we deduce $\|\nabla{\sf p}\|_2=0$, namely,
$$
\nabla p=\nabla\hat{p}\ \ \mbox{in $\Omega$}\,.
$$
From the latter, \eqref{2.71} and \eqref{2.64} we thus infer
\be
\|\nabla p\|_{W^{1,2}(L^2)}\le c\,\mathscr D\,.
\eeq{2.73}
Finally, possibly adding to $p$ a function of time, we observe that
\cite[II.9.1(i)]{Gab}
$$
\|p\|_6\le c \,\|\nabla p\|_2
$$
which, in turn, by \eqref{2.73} and the embedding $W^{1,2}(D^{1,2})\subset L^\infty(D^{1,2})$, proves $p\in \mathscr P$ and the corresponding estimate in \eqref{V29}. The proof of the proposition is completed.
\par\hfill$\square$
\setcounter{equation}{0}
\section{Linear Problem: Existence, Uniqueness and Asymptotic Behavior}
Our next goal is to show that, if $\bfcalf$ also decays ``sufficiently fast" at large spatial distances, then a similar property must hold for the solution $\bfu$ given in \propref{2.1}. In this regard, we notice that, with a suitable choice of the axes, the average of $\bfxi(t)$ over a period can be written as
\be
\lambda\,\bfe_1:=\frac1T\int_0^T\bfxi(t){\rm d}t\,,\ \ \lambda\ge 0\,, 
\eeq{lambda} 
The following result holds.
\Bp Let $\bff$, $\bfxi$  $\bfomega$ and $\bfu$ be as in \propref{2.1}. Suppose in addition that, if $\bfxi(t)\not\equiv \0$, it is $\bfxi(t)=\xi(t)\bfe_1$,  and $\bfomega(t)=\omega(t)\bfe_1$, while no further assumption is imposed on $\bfomega$ if $\bfxi(t)\equiv\0$.
Then, if
$[\!]\bfcalf[\!]_{\infty,2,\lambda}<\infty$, it follows   $[\!]\bfu[\!]_{\infty,1,\lambda}<\infty$ and $p\in L^\infty(L^r)$ for all $r\in (\frac32,6]$. Moreover, there are $C=C(\Omega,T,V_0)$ and $C_1=C_1(\Omega,T,V_0,r)$ such that, setting
\be
\cald:=[\!]\bfcalf[\!]_{\infty,2,\lambda}+\|\bff\|_{W^{1,2}(L^2)}+\|\bfcalf\|_{W^{1,2}(L^2)}+\|\bfxi\|_{W^{2,2}(0,T)}+\|\bfomega\|_{W^{2,2}(0,T)}
\eeq{cald}
we have
\be
[\!]\bfu[\!]_{\infty,1,\lambda}\le C\,\cald\,;\ \ \|p\|_{L^\infty(L^r)}\le C_1\,\cald\,.
\eeq{up}
\EP{2.2}
{\em Proof.} Let $\chi=\chi(s)$ be the ``cut-off" function introduced in \propref{2.1}, and set $\chi_{\bar{R}}(x)=\chi(|x|/\bar{R})$, where $\bar{R}>2R_*$. Further,  let $\bfz$ be a solution to problem 
\eqref{Bog} with $\cala\equiv \{x\in \real^3:\ \half\bar{R}<|x|<\bar{R}\}$, and $f\equiv-\nabla\chi_{\bar{R}}\cdot\bfu$. Since $\int_{\cala}f=0$,  \lemmref{1.1_1} guarantees the existence of such a $\bfz$ with the properties stated there. Setting
\be
\bfw:=\chi_{\bar{R}}\,\bfu+\bfz\,,\ \ {\sf p}:=\chi_{\bar{R}}\,p\,, \ \ \bfcalh=\chi_{\bar{R}}\bfcalf
\eeq{CTM}
from \eqref{4.7} we deduce that $(\bfw,{\sf p})$ is a $T$-periodic solution to the following problem
\be\left.\ba{ll}\medskip
{\partial}_t\bfw-\bfV\cdot\nabla\bfw+\bfomega\times\bfw=\Delta\bfw-\nabla {\sf p}+\Div\bfcalh+\bfg\\
\Div\bfw=0\ea\right\}\ \ \mbox{in ${\real^3}\times (0,T)$}
\,,
\eeq{1}
where
$$
\bfg:=-\bfz_t+\bfV\cdot\nabla\bfz-\bfomega\times\bfz+\Delta\bfz-2\nabla\chi_{\bar{R}}\cdot\nabla\bfu+p\,\nabla\chi_{\bar{R}}-\bfxi\cdot\nabla\chi_{\bar{R}}\,\bfu-\bfu\,\Delta\bfchi_{\bar{R}}-\bfcalf\cdot\nabla\bfchi_{\bar{R}}\,.
$$
Extending $\bfz$ to 0 outside its (spatial) support, we obtain that $\bfg$ is of bounded support as well. From \lemmref{1.1_1},  \propref{2.1} and the assumption on $\bff$  we readily show
\be\ba{ll}\medskip
\Sup{t\ge 0}\|\bfg(t)\|_2\le c\,(\|\bff\|_{W^{1,2}(L^2)}+\|\bfcalf\|_{W^{1,2}(L^2)}+\|\bfxi\|_{W^{2,2}(0,T)}+\|\bfomega\|_{W^{2,2}(0,T)})\,,\\
\Div\bfcalh(t)\in L^\infty(L^2)\,,
\ea
\eeq{g}
where  $c$  depends, at most, on $\Omega$, $T$ and $V_0$. 
Define
$$
\bfA=\bfA(t):=\left(\ba{ccc}\medskip 0&-\omega_3(t) &\omega_2(t)\\ \medskip
\omega_3(t)& 0 &-\omega_1(t)\\
0&\omega_1(t) &0\ea\right)\,,
$$
and let $\bfQ=\bfQ(t)$ be the family of orthogonal transformations that are solutions to the following initial-value problem
$$
\dot{\bfQ}={\bfQ}\cdot\bfA\,,\ \ \bfQ(0)=\bfI\,,
$$ 
with $\bfI$  the identity matrix. It is readily checked that, for any $\bfa\in\real^3$,
\be
\bfQ^\top\cdot\dot{\bfQ}\cdot\bfa=-\dot{\bfQ}^\top\cdot\bfQ\cdot\bfa=\bfomega\times\bfa\,.
\eeq{2.78}
Moreover, in the case when $\omega_2\equiv\omega_3\equiv0$, we also have
\be
\bfQ(t)\cdot\bfe_1=\bfe_1\,,\ \ \mbox{for all $t\ge 0$.}
\eeq{Qq}
We next introduce the new variable $\bfy$ defined by 
\be
\bfy=\bfQ(t)\cdot\bfx+\bfx_0(t)
\eeq{2}
where
\be
\bfx_0(t):=\int_0^t(\bfxi(s)-\lambda{\bfe_1})\,{\rm d}s\,,
\eeq{3}
Since $(1/T)\int_0^T\big(\bfxi(t)-\lambda\,{\bfe_1}\big){\rm d}t=\0$, by \cite[Proposition 1]{GARMA}  there is $M=M(T,V_0)>0$ such that 
\be
\sup_{t\ge 0}|\bfx_0(t)|\le M\,
\eeq{GD0}
which, by \eqref{2} implies, in particular,
\be
|\bfx|-M\le |\bfy|\le |\bfx|+M\,.
\eeq{GD}
Set
\be\ba{ll}\medskip
\bfv(\bfy,t)=\bfQ(t)\cdot\bfw(\bfQ^\top(t)\cdot(\bfy-\bfx_0(t)),t),\\ \medskip
{\sf p}(\bfy,t)=p(\bfQ^\top(t)\cdot(\bfy-\bfx_0(t)),t), 
\\  \medskip
\bfh(\bfy,t)=\bfQ(t)\cdot\bfg(\bfQ^\top(t)\cdot(\bfy-\bfx_0(t)),t)\\
\bfcalg(\bfy,t)=\bfQ(t)\cdot\bfcalh(\bfQ^\top(t)\cdot(\bfy-\bfx_0(t)),t)\cdot\bfQ^\top(t)\,.
\ea
\eeq{4}
By a straightforward calculation we show
$$\ba{rl}\medskip
\pde{\bfv}t=&
\dot{\bfQ}\cdot\bfw+\bfQ\cdot\left\{[\dot{\bfQ}^\top\cdot(\bfy-\bfx_0)-\bfQ^\top\cdot\dot{\bfx}_0]\cdot\nabla\bfw+\pde{\bfw}t\right\}
\\
=&\bfQ\cdot\left\{\bfQ^\top\cdot\dot{\bfQ}\cdot\bfw+[\dot{\bfQ}^\top\cdot(\bfy-\bfx_0)-\bfQ^\top\cdot\dot{\bfx}_0]\cdot\nabla\bfw+\pde{\bfw}t\right\}\,,
\ea
$$
which, in turn, in view of \eqref{2.78}, \eqref{2} and \eqref{3}, provides
\be
\pde{\bfv}t=\bfQ\cdot\left\{\bfomega\times\bfw-[\bfomega\times\bfx+\bfQ^\top\cdot(\bfxi-\lambda\bfe_1)]\cdot\nabla\bfw+\pde{\bfw}t\right\}\,.
\eeq{2.85}
Now, suppose $\bfxi\not\equiv\0$. Then, by assumption and \eqref{Qq} we obtain
\be 
\bfQ^\top\cdot(\bfxi-\lambda\bfe_1)=\bfxi-\lambda\bfe_1\,.
\eeq{2.86}
If, on the other hand, $\bfxi\equiv\0$, then, obviously, $\bfQ^\top\cdot(\bfxi-\lambda\bfe_1)=\bfxi-\lambda\bfe_1=\0$. Thus, in either case, from \eqref{2.85} we deduce
$$
\pde{\bfv}t=\bfQ\cdot\left\{\bfomega\times\bfw-[\bfomega\times\bfx+(\bfxi-\lambda\bfe_1)]\cdot\nabla\bfw+\pde{\bfw}t\right\}\,,
$$ 
and so, by \eqref{1}$_1$,
\be
\pde{\bfv}t=\bfQ\cdot\left(\lambda\,\partial_1\bfw+\Delta\bfw-\nabla p+\Div\bfcalh+\bfg
\right)\,.
\eeq{2.87} 
By a straightforward calculation, we show
$$
\Delta_y\bfv=\bfQ\cdot\Delta_x\bfw\,,\ \nabla_y{\sf p}=\bfQ\cdot\nabla_xp\,,\ \ \Div_y\bfcalg=\bfQ\cdot\Div_x\bfcalh\,,\ \ \Div_y\bfv=\Div_x\bfw\,, 
$$
and also, by assumption and \eqref{Qq},
$$
\partial_1\bfv=\bfQ\cdot\partial_1\bfw\,.
$$ 
As a result, from \eqref{2.87} we infer that $(\bfv,{\sf p})$ is a solution to the following Cauchy problem
\be\ba{cc}\medskip\left.\ba{ll}\medskip
\bfv_t-\lambda\,\partial_1\bfv=\Delta\bfv-\nabla {\sf p}+\Div\bfcalg+\bfh\\
\Div\bfv=0\ea\right\}\ \ \mbox{in ${\real^3}\times (0,\infty)$}
\,,\\
\bfv(x,0)=\bfw(x,0)\,.
\ea
\eeq{5}
In order to obtain the desired spatial decay property, we would like to apply \lemmref{VIII.4.4} to \eqref{5}. To this end, we begin to observe that, by \eqref{esse} and \eqref{GD} it follows that
\be\ba{rl}\medskip
(1+|x|)(1+2\lambda\,s(x))&\le (1+|y|+M)\big(1+2\lambda\,s(y)+2\lambda\,(M+x_{01}(t))\big)
\\
&\le c\,(1+|y|)\,\big(1+2\lambda\,s(y)\big)\,,
\ea
\eeq{LJ}
and, likewise,
\be
(1+|y|)(1+2\lambda\,s(y))\le c\, (1+|x|)\big(1+2\lambda\,s(x)\big)
\,.
\eeq{VFM}
We next look for a solution to \eqref{5} of the form $(\bfv_1+\bfv_2, {\sf p})$ where
\be\ba{cc}\medskip\left.\ba{ll}\medskip
(\bfv_1)_t-\lambda\,\partial_1\bfv_1=\Delta\bfv_1-\nabla {\sf p}_1+\Div\bfcalg+\bfh\\
\Div\bfv_1=0\ea\right\}\ \ \mbox{in ${\real^3}\times (0,\infty)$}
\,,\\
\bfv_1(x,0)=\0\,,
\ea
\eeq{6}
and
\be\ba{cc}\medskip\left.\ba{ll}\medskip
(\bfv_2)_t-\lambda\,\partial_1\bfv_2=\Delta\bfv_2\\
\Div\bfv_2=0\ea\right\}\ \ \mbox{in ${\real^3}\times (0,\infty)$}
\,,\\
\bfv(x,0)=\bfw(x,0)\,.
\ea
\eeq{7}
From $\eqref{4}_{4}$,  \eqref{VFM} and \eqref{CTM}$_3$ we infer
\be 
[\!]\bfcalh[\!]_{\infty,2,\lambda}\le C\,[\!]\bfcalf[\!]_{\infty,2,\lambda}
\eeq{1.45}
Moreover, by \eqref{4}$_3$ and \eqref{g}, it follows that 
\be
\Sup{t\ge 0}\|\bfh(t)\|_2\le c\,(\|\bff\|_{W^{1,2}(L^2)}+\|\bfcalf\|_{W^{1,2}(L^2)}+\|\bfxi\|_{W^{2,2}(0,T)}+\|\bfomega\|_{W^{2,2}(0,T)})\,.
\eeq{1.46}
As a result, from \lemmref{VIII.4.4} we conclude that \eqref{7}
has one and only one solution such that for all $T>0$,
$$\ba{ll}\medskip
\bfv_1\in L^2(0,T;W^{2,2})\,, \ (\bfv_1)_t\in L^2(0,T;L^2)\,;\ \ \nabla {\sf p}\in L^{2}(0,T;L^2)\,;\\  
{} [\!] \bfv_1 [\!]_{\infty,1,\lambda}<\infty\,,\ {\sf p}_1\in L^\infty(L^r)\,,\, \mbox{for all $r\in (\frac32,\infty)$}
\ea
$$
satisfying, in addition, the inequality
\be \ba{rl}\medskip
[\!] \bfv_1 [\!]_{\infty,1,\lambda}\le C\,\cald\,,\ \ \  
\| {\sf p} \|_{L^\infty(L^r)} \le C_1 \,\cald\,, 
\ea
\eeq{munn2}
with $C$ and $C_1$ as in \eqref{up}.
Concerning \eqref{7}, 
since by \propref{2.1}, \eqref{CTM} and \eqref{Bog} it is $\bfw(x,0)\in L^6(\real^3)$ and $\Div\bfw(x,0)=0$, it follows that there exists a (unique) solution $\bfv_2 $ such that (see, e.g., \cite[Theorem VIII.4.3]{Gab})
\be\ba{ll}\medskip
\bfv_2,\partial_t\bfv_2\, D^2\bfv_2\in L^r([\varepsilon,\tau]\times\real^3)\,,\ \ \mbox{all $\varepsilon\in (0,\tau)$, $\tau>0$, and $r\in[6,\infty]$}\,,\\ 
\|\bfv_2(t)\|_\infty\le C_2\,t^{-\frac14}\|\bfw(0)\|_6\,,\ \ \Sup{t\in(0,\infty)}\|\bfv_2(t)\|_6\le C_2\,\|\bfw(0)\|_6\,.
\ea 
\eeq{15}
In view of the regularity properties of $\bfu$ (and hence of $\bfw$) and those in (\ref{estimunn1}), \eqref{15} for $\bfv_i$, $i=1,2$, respectively, we may use the results proved in \cite[Lemma VIII.4.2]{Gab} to guarantee, by uniqueness,  that $\bfv_1+\bfv_2\equiv\bfv$ and ${\sf p}_1\equiv {\sf p}$. Thus, in particular, from \eqref{munn2} we get
\be\ba{ll}\medskip
\|p\|_{L^\infty(L^r(\Omega^{\bar{R}})}\le c\,\left({}[\!]\bfcalf[\!]_{\infty,2,\lambda} +\|\bff\|_{W^{1,2}(L^2)}+\|\bfcalf\|_{W^{1,2}(L^2)}\right.\\ \hspace*{3.5cm}\left.+\|\bfxi\|_{W^{2,2}(0,T)}+\|\bfomega\|_{W^{2,2}(0,T)}\right)\,,\ \ \mbox{for all $r\in (\frac32,\infty)$}\,,
\ea
\eeq{press}
where $c=c(\Omega,T,V_0,r)$.
Further, 
due to the $T$-periodicity of $\bfw$ and \eqref{4}$_1$, for  arbitrary positive integer $n$ and $t\in[0,T]$ we obtain 
\be\ba{rl}\medskip
|\bfw(x,t)|(1+|x|)(1+&2\lambda\,s(x))=|\bfv(y,t+nT)|(1+|x|)(1+2\lambda\,s(x))\\ &\le \big(|\bfv_1(y,t+nT)|+|\bfv_2(y,t+nT)|\big)(1+|x|)(1+2\lambda\,s(x)).\ea
\eeq{17}
Employing  \eqref{LJ}, \eqref{munn2} and \eqref{15}$_2$ in this inequality we get
$$\ba{ll}\medskip
|\bfw(x,t)|(1+|x|)(1+2\lambda\,s(x))
\le c\,\big[(1+|x|)(1+2\lambda\,s(x)) (t+nT)^{-\frac14}\|\bfw(0)\|_6\\
\hspace*{1.5cm} +[\!]\bfcalf[\!]_{\infty,2,\lambda} +\|\bff\|_{W^{1,2}(L^2)}+\|\bfcalf\|_{W^{1,2}(L^2)}+\|\bfxi\|_{W^{2,2}(0,T)}+\|\bfomega\|_{W^{2,2}(0,T)}\big]\,\ea
$$ 
so that,  by letting $n\to\infty$ and recalling that, uniformly in $t\ge 0$,  $\bfu(x,t)\equiv \bfw(x,t)$  for $|x|$ sufficiently large ($>\bar{R}$)  we deduce
\be
[\!]\bfu[\!]_{\infty,1,\lambda,\Omega^{\bar{R}}}\le c\,\big([\!]\bfcalf[\!]_{\infty,2,\lambda} +\|\bff\|_{W^{1,2}(L^2)}+\|\bfcalf\|_{W^{1,2}(L^2)}+\|\bfxi\|_{W^{2,2}(0,T)}+\|\bfomega\|_{W^{2,2}(0,T)}\big)\,. 
\eeq{vfn}
Moreover, by \propref{2.1}  we have 
\be
\|\bfu\|_{L^\infty(L^\infty)}+\|p\|_{L^\infty(L^6)}\le c\,\big(\|\bff\|_{W^{1,2}(L^2)}+\|\bfcalf\|_{W^{1,2}(L^2)}+\|\bfxi\|_{W^{2,2}(0,T)}+\|\bfomega\|_{W^{2,2}(0,T)}\big)
\eeq{vfm}
and 
the desired result then follows from  \eqref{press}, \eqref{vfn} and \eqref{vfm}.
\par\hfill$\square$\par
Combining the  results of \propref{2.1} and \propref{2.2} we can prove the   main achievement of this section. Precisely, let
$$\ba{ll}\medskip
{\mathscr U}_\lambda:=\{\bfu\in\mathscr U: \ \|\bfu\|_{\mathscr U_\lambda}:=\|\bfu\|_{\mathscr U}+[\!]\bfu[\!]_{\infty,1,\lambda}<\infty\}\,,\\
{\mathscr P}:=\{p\in\hat{\mathscr P}:\ \|p\|_{{\mathscr P}}:=\|p\|_{\hat{\mathscr P}}+\|p\|_{L^\infty(L^r)}<\infty\,,\ \mbox{for all $r\in(\frac32,6]$}\}\,,
\ea
$$
with $\lambda$ defined in \eqref{lambda}. 
The following theorem holds.
\Bt Let $\bff=\Div\bfcalf$, and $\bff,\bfcalf\in W^{1,2}(L^2)$,   $\bfxi,\bfomega\in W^{2,2}(0,T)$ be given $T$-periodic functions, with $[\!]\bfcalf[\!]_{\infty,2,\lambda}<\infty$.
Suppose in addition that, if $\bfxi(t)\not\equiv \0$, it is   $\bfxi(t)=\xi(t)\bfe_1$,  and $\bfomega(t)=\omega(t)\bfe_1$, while no further assumption is imposed on $\bfomega$ if $\bfxi(t)\equiv\0$. Then,  there exists one and only one $T$-periodic solution $(\bfu,p)$ to \eqref{4.7} with $(\bfu,p)\in \mathscr U_\lambda\times\mathscr P$.
Moreover, the following estimate holds
$$
\|\bfu\|_{\mathscr U_\lambda}\le C\,\cald\,,\ \ \|p\|_{\mathscr P}
\le C_1 \,\cald
$$
where $\cald$ is defined in \eqref{cald},  $C=C(\Omega, T,V_0)$, $C_1=C_1(\Omega,T,V_0,r)$, and $V_0$ is given in \eqref{V_0}.
\ET{1.1}
{\em Proof.} In view of \propref{2.1} and \propref{2.2}, it remains to show the uniqueness property, namely, that, under the given assumption on $\bfV$, the problem
\be\ba{cc}\smallskip\left.\ba{ll}\medskip
\bfu_t-\bfV(t)\cdot\nabla\bfu+\bfomega(t)\times\bfu=\Delta\bfu-\nabla {p}\\
\Div\bfu=0\ea\right\}\ \ \mbox{in $\Omega\times (0,T)$}\\
\bfu(x,t)=\0\,,\ \ (x,t)\in \partial\Omega\times [0,T]\,,
\ea
\eeq{Uni}
has only the null solution in the class $(\bfu,p)\in\mathscr U_\lambda\times\mathscr P$. To this end, we recall \cite[Lemma 3]{GaSiS} \cite[Lemma II.6.4]{Gab} that there exists a ``cut-off" function, $\psi_R\in C_0^\infty(\real^3)$, $R\in(0,\infty)$,  with the following properties
\begin{itemize}
  \item [(ii)] $\psi_R(x)\in [0,1]$, for all $x\in\real^3$ and $R>0$;
  \item [(ii)] $\Lim{R\to\infty}\psi_R(x)=1$ for all $x\in\real$;
  \item [(iii)] $\nabla\psi_R(x)\cdot(\bfomega\times\bfx)=0$, for all $x\in\real$ and $R>0$;
  \item [(iv)] $\supp(\psi_R)\subset \{x\in\real^3:2R<|x|<(2R)^2\}$;
  \item [(v)] $\|\bfu\,\nabla\psi_R\|_2\le c\,\|\nabla\bfu\|_{2,\Omega^{\frac{R}{\sqrt2}}}$,\, with $c$ independent of $R$;
  \item [(vi)] $|x|^{-2}\partial_1\psi_R\in L^1(\Omega)$,\, and $\Lim{R\to\infty}\Int{\Omega}{}\Frac{|\partial_1\psi_R|}{|x|^2}=0$\,.
\end{itemize}
Let us dot-multiply both sides of \eqref{Uni}$_1$ by $\psi_R\bfu$ and integrate by parts over $\Omega\times[0,T]$. Using $T$-periodicity, (iii)  and \eqref{Uni}$_{2,3}$, we thus obtain
\be\ba{rl}\medskip
\Int0T\|\sqrt{\psi_R}\,\nabla\bfu\|_{2}^2&
=\Int0T\Int{\Omega}{}\left(-\nabla\psi_R\cdot\nabla\bfu\cdot\bfu+\half\xi\partial_1\psi_R|\bfu|^2+p\nabla\psi_R\cdot\bfu\right)\\
&
:= I_1(R)+I_2(R)+I_3(R)\,.
\ea
\eeq{3.28}
By Schwarz inequality,  (iv) and (v), we show 
\be
|I_1(R)|\le \|\bfu\nabla\psi_R\|_{2,\Omega^{\frac R{\sqrt2}}}\|\nabla\bfu\|_{2,\Omega^{\frac R{\sqrt2}}}\le c\,\|\nabla\bfu\|_{2,\Omega^{\frac R{\sqrt2}}}^2\,,
\eeq{3.30}
and, likewise,
\be
|I_3(R)|\le c\,\|p\|_{2,\Omega^{\frac R{\sqrt2}}}\,\|\nabla\bfu\|_{2,\Omega^{\frac R{\sqrt2}}}\,.
\eeq{3.31}
Since $\bfu(x,t)=O(|x|^{-1})$ uniformly in time, we also get
\be
|I_2(R)|\le c\,\int_0^T\int_\Omega\frac{|\partial_1\psi_R|}{|x|^2}\,.
\eeq{3.32}
Thus, letting $R\to\infty$ in \eqref{3.28}, and  taking into account \eqref{3.30}--\eqref{3.32}, (i), (vi) and  that $(\bfu,p)\in \mathscr U_\lambda\times\mathscr P$, by the dominated convergence theorem we infer
$$
\int_0^T\|\nabla\bfu(t)\|^2_2=0\,,
$$
which, in view of \eqref{Uni}$_3$ and the summability properties of $p$, furnishes $\bfu\equiv p\equiv 0$.
\par\hfill$\square$\par
\setcounter{equation}{0}
\section{On the Unique Solvability of the Nonlinear Problem}
In this section we shall prove existence and uniqueness of $T$-periodic solutions to the  full nonlinear problem \eqref{0.1}, provided the magnitude of the data is suitably restricted. 
This goal will be reached by combining  \theoref{1.1} with a contraction mapping argument. 
\par 
Set
$$
\mathscr F_\lambda:=\{\bfcalf:\Omega\times[0,T]\mapsto \real^{3\times3}\,:\, \bfcalf\,,\,\Div\bfcalf\in W^{1,2}(L^2),\, {}[\!]\bfcalf[\!]_{\infty,2,\lambda}<\infty\}\,, 
$$
with $\lambda$ defined in \eqref{lambda}, and define
$$
\mathscr D_\lambda:=\{ (\bfcalf,\bfxi,\bfomega)\,\ \mbox{$T$-periodic}:\, \bfcalf\in \mathscr F_\lambda;\,\ \bfxi,\bfomega\in W^{2,2}(0,T)\}
$$
endowed with the norm
$$
\|(\bfcalf,\bfxi,\bfomega)\|_{\mathscr D_\lambda}:=\|\Div\bfcalf\|_{W^{1,2}(\Omega)}+\|\bfcalf\|_{W^{1,2}(\Omega)}+[\!]\bfcalf[\!]_{\infty,2,\lambda}+\|\bfxi\|_{W^{2,2}(0,T)}+\|\bfomega\|_{W^{2,2}(0,T)}\,.
$$
We need the following preliminary result.
\Bl Let $\textbf{\textsf{u}},\textbf{\textsf{w}}\in \mathscr U_\lambda$. Then $\textbf{\textsf{u}}\otimes\textbf{\textsf{w}}\in\mathscr F_\lambda$ and
\be
\|\Div(\textbf{\textsf{u}}\otimes \textbf{\textsf{w}})\|_{W^{1,2}(L^2)}+\|\textbf{\textsf{u}}\otimes \textbf{\textsf{w}}\|_{W^{1,2}(L^2)}+[\!]\textbf{\textsf{u}}\otimes\textbf{\textsf{w}}[\!]_{\infty,2,\lambda}\le c\|\textbf{\textsf{u}}\|_{\mathscr U_\lambda}\|\textbf{\textsf{w}}\|_{\mathscr U_\lambda}.
\eeq{esti}
\EL{3.1}
{\em Proof.} Obviously,
\be
[\!]\textbf{\textsf{u}}\otimes\textbf{\textsf{w}}[\!]_{\infty,2,\lambda}\le \|\textbf{\textsf{u}}\|_{\mathscr U_\lambda}\|\textbf{\textsf{w}}\|_{\mathscr U_\lambda}\,.
\eeq{conti0}
We next observe that, since $\Div\textbf{\textsf{u}}=0$, it follows 
\be
\Div(\textbf{\textsf{u}}\otimes \textbf{\textsf{w}})=\textbf{\textsf{u}}\cdot\nabla\textbf{\textsf{w}}.
\eeq{dive} 
Now, clearly,
\be
\|\textbf{\textsf{u}}\cdot\nabla\textbf{\textsf{w}}\|_{L^{2}(L^2)}\le [\!]\textbf{\textsf{u}}[\!]_{\infty,1,\lambda}\,\|\nabla\textbf{\textsf{w}}\|_{L^\infty(L^2)}\le \|\textbf{\textsf{u}}\|_{\mathscr U_\lambda}\|\textbf{\textsf{w}}\|_{\mathscr U_\lambda}\,,
\eeq{conti}
and
\be
\|\textbf{\textsf{u}}\otimes\textbf{\textsf{w}}\|_{L^{2}(L^2)}\le c\,[\!]\textbf{\textsf{u}}[\!]_{\infty,1,\lambda}\,[\!]\textbf{\textsf{w}}[\!]_{\infty,1,\lambda}\le c\,\|\textbf{\textsf{u}}\|_{\mathscr U_\lambda}\|\textbf{\textsf{w}}\|_{\mathscr U_\lambda}\,.
\eeq{conti1}
Moreover, by using the inequality $\|\nabla w\|_{3}\le c\,\|\nabla w\|_2^{\frac12}\|D^2w\|_2^{\frac12}$ (see \cite[Theorem 2.1]{CrMa}) along with H\"older inequality, we get
\be\ba{rl}\medskip
\|\textbf{\textsf{u}}_t\cdot\nabla\textbf{\textsf{w}}\|_{L^{2}(L^2)}\!+\!\|\textbf{\textsf{u}}\cdot\nabla\textbf{\textsf{w}}_t&\!\|_{L^{2}(L^2)}\le\|\textbf{\textsf{u}}_t\|_{L^{\infty}(L^6)} \|\nabla\textbf{\textsf{w}}\|_{L^2(L^3)}+[\!]\textbf{\textsf{u}}[\!]_{\infty,1,\lambda}\, \|\nabla\textbf{\textsf{w}}_t\|_{L^2(L^2)}\\ \medskip
&\!\!\! \le c\,\big(\|\textbf{\textsf{u}}\|_{\mathscr U_\lambda}(\|\nabla\textbf{\textsf{w}}\|_{L^\infty(L^{2})}^{\frac12}\|D^2\textbf{\textsf{w}}\|_{L^2(L^{2})}^{\frac12})+\|\textbf{\textsf{u}}\|_{\mathscr U_\lambda}\|\textbf{\textsf{w}}\|_{\mathscr U_\lambda}\big)\\
&\!\!\! \le c\,\|\textbf{\textsf{u}}\|_{\mathscr U_\lambda}\|\textbf{\textsf{w}}\|_{\mathscr U_\lambda}\,.
\ea
\eeq{conti2}
Thus, combining \eqref{conti2} and \eqref{conti} and \eqref{dive} we conclude
\be
\|\Div(\textbf{\textsf{u}}\otimes \textbf{\textsf{w}})\|_{W^{1,2}(L^2)}\le c\|\textbf{\textsf{u}}\|_{\mathscr U_\lambda}\|\textbf{\textsf{w}}\|_{\mathscr U_\lambda}.
\eeq{conti3}
By Hardy inequality \cite[Theorem II.6.1]{Gab}, we infer
$$\ba{rl}\medskip
\|\textbf{\textsf{u}}_t\otimes\textbf{\textsf{w}}\|_{L^{2}(L^2)}+\|\textbf{\textsf{u}}\otimes\textbf{\textsf{w}}_t\|_{L^{2}(L^2)}&\!\!\!\le [\!]\textbf{\textsf{w}}[\!]_{\infty,1,\lambda}\|\textbf{\textsf{u}}_t/|x|\|_{L^2(L^2)}+[\!]\textbf{\textsf{u}}[\!]_{\infty,1,\lambda}\|\textbf{\textsf{w}}_t/|x|\|_{L^2(L^2)}\\ \medskip
&\!\!\!\le c\,\left([\!]\textbf{\textsf{w}}[\!]_{\infty,1,\lambda}\|\nabla\textbf{\textsf{u}}_t\|_{L^2(L^2)}+[\!]\textbf{\textsf{u}}[\!]_{\infty,1,\lambda}\|\nabla\textbf{\textsf{w}}_t\|_{L^2(L^2)}\right)\\
&\!\!\!\le c\,\|\textbf{\textsf{u}}\|_{\mathscr U_\lambda}\|\textbf{\textsf{w}}\|_{\mathscr U_\lambda}\,.
\ea
$$
The latter and \eqref{conti1} furnish
\be
\|\textbf{\textsf{u}}\otimes\textbf{\textsf{w}}\|_{W^{1,2}(L^2)}\le c\,\|\textbf{\textsf{u}}\|_{\mathscr U_\lambda}\|\textbf{\textsf{w}}\|_{\mathscr U_\lambda}\,.
\eeq{conti4}
The lemma follows from \eqref{conti0},  \eqref{conti3} and \eqref{conti4}.
\par\hfill$\square$\par
We are now in a position to prove the main result of this paper.
\Bt
Let $(\bfcalb, \bfxi,\bfomega)\in\mathscr D_\lambda$.
Suppose in addition that, if $\bfxi(t)\not\equiv \0$, it is   $\bfxi(t)=\xi(t)\bfe_1$,  and $\bfomega(t)=\omega(t)\bfe_1$, while no further assumption is imposed on $\bfomega$ if $\bfxi(t)\equiv\0$.
Then, there exists  $\varepsilon_0=\varepsilon_0(\Omega,T,V_0)>0$\footnote{$V_0$ defined in \eqref{V_0}.} such that if 
$$
\|(\bfcalb, \bfxi,\bfomega)\|_{\mathscr D_\lambda}<\varepsilon_0\,,
$$
problem \eqref{0.1} has
one and only one   solution $(\bfu,p)\in \mathscr U_\lambda\times
 \mathscr P$\, with  $\|\bfu\|_{\mathscr U_\lambda}+\|p\|_{\mathscr P}\le c\,\|(\bfcalb, \bfxi,\bfomega)\|_{\mathscr D_\lambda}$, for some $c=c(\Omega,T,V_0,r)$.
\ET{3.1}
{\em Proof.} We want to apply  the contraction mapping theorem to the map 
$$
M:\textbf{\textsf{u}}\in\mathscr U_\lambda\mapsto \bfu\in\mathscr U_\lambda\,,
$$
with $\bfu$ solving the linear problem
\be\ba{cc}\smallskip\left.\ba{ll}\medskip
\bfu_t-\bfV\cdot\nabla\bfu+\bfomega\times\bfu=\Delta\bfu-\nabla {p}+\textbf{\textsf{u}}\cdot\nabla \textbf{\textsf{u}}+\Div\bfcalb\\
\Div\bfu=0\ea\right\}\ \ \mbox{in $\Omega\times (0,T)$}\\
\bfu(x,t)=\bfV(t)\,,\ \ (x,t)\in \partial\Omega\times [0,T]\,,
\ea
\eeq{lin}
Set 
\be\textbf{\textsf{u}}\cdot\nabla\textbf{\textsf{u}}=\Div(\textbf{\textsf{u}}\otimes \textbf{\textsf{u}}):=\Div\textbf{\textsf{F}}\,,
\eeq{C}
where we used the condition $\Div\textbf{\textsf{u}}=0$. In virtue of \lemmref{3.1},  by assumption, and by the obvious inequality
$$
[\!]\textbf{\textsf{F}}[\!]_{\infty,2,\lambda}\le c_1 [\!]\textbf{\textsf{u}}[\!]_{\infty,1,\lambda}^2\,,\ \ \textbf{\textsf{u}}\in\mathscr U_\lambda\,,
$$
we infer that
$\textbf{\textsf{F}}$,  $\bfcalb$,  $\bfxi$ and $\bfomega$ satisfy the assumptions of \theoref{1.1}. Therefore, by that theorem we conclude that the map $M$ is well defined and, in particular, that
\be
\|\bfu\|_{\mathscr U_\lambda}\le c_2\left(\|\textbf{\textsf{u}}\|_{\mathscr U_\lambda}^2+\|(\bfcalb,\bfxi,\bfomega)\|_{\mathscr D_\lambda}\right)\,,
\eeq{3.3}
with $c_2=c_2(\Omega,T,V_0)$. If we now take \be\|\textbf{\textsf{u}}\|_{\mathscr U_\lambda}<\delta\,,\ \  \delta:=4c_2\|(\bfcalb,\bfxi,\bfomega)\|_{\mathscr D_\lambda}\,,\ \ \|(\bfcalb,\bfxi,\bfomega)\|_{\mathscr D_\lambda}<\frac{1}{16 c_2^2}\,,
\eeq{BD}
from \eqref{3.3} we deduce $\|\bfu\|_{\mathscr U_\lambda}<\half\delta$.
Let $\textbf{\textsf{u}}_i\in\mathscr U_\lambda$ $i=1,2$, and set 
$$
\textbf{\textsf{u}}:=\textbf{\textsf{u}}_1-\textbf{\textsf{u}}_2\,,\ \ \bfu:=M(\textbf{\textsf{u}}_1)-M(\textbf{\textsf{u}}_2)\,.
$$
From \eqref{lin} we then show
\be\ba{cc}\smallskip\left.\ba{ll}\medskip
\bfu_t-\bfV\cdot\nabla\bfu+\bfomega\times\bfu=\Delta\bfu-\nabla {p}+\textbf{\textsf{u}}_1\cdot\nabla \textbf{\textsf{u}}+\textbf{\textsf{u}}\cdot\nabla\textbf{\textsf{u}}_2\\
\Div\bfu=0\ea\right\}\ \ \mbox{in $\Omega\times (0,T)$}\\
\bfu(x,t)=\0\,,\ \ (x,t)\in \partial\Omega\times [0,T]\,.
\ea
\eeq{line}
Proceeding as in the proof of \eqref{3.3} we can show
$$
\|\bfu\|_{\mathscr U_\lambda}\le c_2\,\left(\|\textbf{\textsf{u}}_1\|_{\mathscr U_\lambda}+\|\textbf{\textsf{u}}_2\|_{\mathscr U_\lambda}\right)\|\textbf{\textsf{u}}\|_{\mathscr U_\lambda}\,.
$$
As a result, since $\|\textbf{\textsf{u}}_i\|_{\mathscr U_\lambda}<\delta$, $i=1,2$, from the previous inequality we infer
$$
\|\bfu\|_{\mathscr U_\lambda}< 2c_2\delta\|\textbf{\textsf{u}}\|_{\mathscr U_\lambda}\,.
$$
By \eqref{BD}, we have $2c_2\delta<1/2$, and so, from the last displayed relation we may conclude that $M$ is a contraction, which, along with \eqref{BD}, completes the proof of the theorem.
\par\hfill$\square$\par

\ed